\documentclass[12pt]{article}

\usepackage[]{pdfsync} 
\usepackage[francais,english]{babel} 
\usepackage[]{latexsym} 
\usepackage{amsmath} 
\usepackage{graphicx} 
\usepackage{amsfonts} 
\usepackage{amssymb}
\usepackage[T1]{fontenc}

\newcommand \as{almost surely}  

\font \db = msbm10 at 12 pt
 
 \def\iff{if and only if } 

       \def \d{\delta}     \def \G{\Gamma}
        
\def \e{\varepsilon}  \def \f{\varphi}    \def \rr{\varrho}

\font \db = msbm10 at 12 pt
     \font \bdb = msbm10 at 14 pt  
\def \R{\mbox{\db R}}  \def \H{\mbox{\db H}} \def \E{\mbox{\db E}} \def \P{\mbox{\db P}}

\def \bR{\mbox{\bdb R}}       

\def \Tr{\hbox{\rm Tr}}   

\newcommand \parsn{\par \smallskip \noindent }  \newcommand \parmn{\par \medskip \noindent } 
\newcommand \pars{\par \smallskip }  \newcommand \parm{\par \medskip }

\def \S{\mbox{\db S}}

\def \AA{{\cal A}}           \def \EE{{\cal E}}
    \def \HH{{\cal H}}       
\def \LL{{\cal L}}  \def \M{{\cal M}}       \def \O{{\cal O}}    
  \def \RR{{\cal R}}      \def \TT{{\cal T}}   
  \def \KK{{\cal K}}

\def \p{\partial}            
  \def \lra{\longrightarrow}   \def \ub{\underbar}

              \def \nea{\nearrow}        

   \def \ds{\displaystyle}   \def \ts{\textstyle}       
\def \rt1{\sqrt{-1}\,\,}  \def \1{^{-1}}            \def \2{^{-2}}             \def \5{{\ts \frac{1}{2}}}
    \def\indf{\leavevmode\indent }

\def \parn{\par\noindent }    \def\indf{\leavevmode\indent }  
\def \parsn{\par\smallskip \noindent } \def \parmn{\par\medskip \noindent } 
\def \parm{\par \medskip } \def \pars{\par \smallskip }

\begin{document}

\newtheorem{defi}[subsubsection]{Definition} 
\newtheorem{theo}[subsubsection]{Theorem} 
\newtheorem{prop}[subsubsection]{Proposition}  
\newtheorem{cor}[subsubsection]{Corollary} 
\newtheorem{lem}[subsubsection]{Lemma} 
\newtheorem{rem}[subsubsection]{Remark} 
\newtheorem{subs}[subsubsection]{} 

\newtheorem{Cor}[subsection]{Corollary} 
\newtheorem{Lem}[subsection]{Lemma} 
\newtheorem{Theo}[subsection]{Theorem} 
\newtheorem{Prop}[subsection]{Proposition} 
\newtheorem{Rem}[subsection]{Remark} 

\newcommand \beq{\begin{equation}} \newcommand \eeq{\end{equation}} 

\newcommand \bthe{\begin{theo}}  \newcommand \ethe{\end{theo}}   
\newcommand \bpro{\begin{prop}}  \newcommand \epro{\end{prop}}   
\newcommand \bcor{\begin{cor}}    \newcommand \ecor{\end{cor}}     
\newcommand \blem{\begin{lem}}   \newcommand \elem{\end{lem}}   
\newcommand \brem{\begin{rem}}  \newcommand \erem{\end{rem}} 
\newcommand \bdefi{\begin{defi}}   \newcommand \edefi{\end{defi}}  
\newcommand \bsub{\begin{subs} \  }   \newcommand \esub{\end{subs}}  

\newcommand \bCor{\begin{Cor}}    \newcommand \eCor{\end{Cor}}     
\newcommand \bLem{\begin{Lem}}   \newcommand \eLem{\end{Lem}}   
\newcommand \bThe{\begin{Theo}}  \newcommand \eThe{\end{Theo}}   
\newcommand \bPro{\begin{Prop}}  \newcommand \ePro{\end{Prop}}   
\newcommand \bRem{\begin{Rem}}  \newcommand \eRem{\end{Rem}}

\title{\bf Curvature Diffusions in General Relativity}
\author{Jacques FRANCHI \quad and \quad Yves LE JAN}
\date {March 2011}
\maketitle

\begin{abstract}
We define and study on Lorentz manifolds a family of covariant diffusions in which the quadratic variation  is locally determined by the curvature. This allows the interpretation of the diffusion effect on a particle by its interaction with the ambient space-time. We will focus on the case of warped products, especially Robertson-Walker manifolds, and analyse their asymptotic behaviour in the case of Einstein-de Sitter-like manifolds.     
\end{abstract}

\centerline{\bf  Contents} \parsn 
\ref{sec.I}) \ Introduction  \hfill  page \pageref{sec.I} \parn 
\ref{sec.IntCcurvRelDiff})  \  Canonical vector fields and curvature  \hfill  page \pageref{sec.IntCcurvRelDiff} \par 
    \ref{sec.IsomL^2Rso1d}) \  Isomorphism between $\bigwedge^2\R^{1,d}$ and $\,so(1,d)$   \hfill  page \pageref{sec.IsomL^2Rso1d}   \par     
    \ref{sec.FramebGM}) \  Frame bundle $\,G(\M)$ over $(\M,g)$    \hfill  page 
      \pageref{sec.FramebGM} \par 
    \ref{sec.ExprLocCoord}) \  Expressions in local coordinates\hfill  page \pageref{sec.ExprLocCoord} \par 
\ref{sec.perffluid}) \  Case of a perfect fluid  \hfill  page \pageref{sec.perffluid} \parn 
\ref{sec.GCCRD}) \  Covariant $\,\Xi$-relativistic diffusions \hfill  page \pageref{sec.GCCRD}  \par 
    \ref{sec.reldiff}) \   The basic relativistic diffusion   \hfill  page \pageref{sec.reldiff}\par    
\ref{sec.XiGCCRD}) \  Construction of the $\,\Xi$-diffusion \hfill  page \pageref{sec.XiGCCRD} \par    
   \ref{sec.RLRD}) \  The $R$-diffusion \hfill  page \pageref{sec.RLRD} \par    
   \ref{sec.RicciRelD}) \  The $\,\EE$-diffusion  \hfill  page \pageref{sec.RicciRelD} \parn
\ref{sec.CWP}) \  Warped (or skew) products  \hfill  page \pageref{sec.CWP} \parn  
\ref{sec.exRW}) \  Example of Robertson-Walker manifolds  \hfill  page \pageref{sec.exRW} \par  
    \ref{sec.RWR1}) \  $\Xi$-relativistic diffusions in an Einstein-de Sitter-like manifold  \hfill  page \pageref{sec.RWR1} \par 
\ref{sec.SasSt}) \   Asymptotic behavior of the $\,R$-diffusion in an E.-d.S.-like manifold  \hfill  page \pageref{sec.SasSt}  \par  
\ref{sec.BehavtCED}) \  Asymptotic energy of the $\,\EE$-diffusion  in an E.-d.S.-like manifold  \hfill  page \pageref{sec.BehavtCED}    \parn  
     \ref{sec.SRD}) \  Sectional relativistic diffusion  \hfill  page \pageref{sec.SRD} \par 
\ref{sec.IntgenGM}) \  Intrinsic relativistic generators on $\,G(\M)$  \hfill  page \pageref{sec.IntgenGM} \par 
    \ref{sec.PosCondCurv}) \  Sign condition on timelike sectional curvatures  \hfill  page \pageref{sec.PosCondCurv} \par
    \ref{sec.RWR12}) \  Sectional diffusion in an Einstein-de Sitter-like manifold  \hfill  page \pageref{sec.RWR12}      \parn  
\ref{sec.S}) \  References  \hfill  page \pageref{sec.S}  

\section{Introduction} \label{sec.I} \indf 
   It is known since Dudley's pioneer's work [Du] that a relativistic diffusion, i.e. a Lorentz-covariant Markov diffusion process, cannot exist on the base space, even in the Minkowski framework of special relativity, but possibly makes sense at the level of the tangent bundle. In this spirit, the general case of a Lorentz manifold was first investigated in [F-LJ], where a general relativistic diffusion was introduced. The quadratic variation of this diffusion is constant 
and does not vanish in the vacuum. In this article, we investigate diffusions whose quadratic variation is locally determined by the curvature of the space, and then vanishes in empty (or at least flat) regions.  \par
   
   The relativistic diffusion considered in [F-LJ] lives on the pseudo-unit tangent bundle $T^1\M$ of the given generic Lorentz manifold $(\M,g)$. As is recalled in Section \ref{sec.reldiff} below, it can be obtained by superposing, to the geodesic flow of $T^1\M$, random fluctuations of the velocity that are given by hyperbolic Brownian motion, if one identifies the tangent space $\,T^1_\xi\M\,$ with the hyperbolic space $\,\H^d$ (at point $\,\xi\in\M$, by means of the pseudo-metric $\,g$). These Brownian fluctuations of the velocity can be defined by the vertical Dirichlet form \  ${\ds \int_{T^1\M} \Big|\nabla_{\dot\xi}^vF(\xi,\dot\xi) \Big|^2\, \mu(d\xi,d\dot\xi)\,}$, considered with respect to the Liouville measure $\,\mu\,$. \pars  
   
   The Dirichlet forms we investigate in this article depend only on the local geometry of $(\M,g)$, e.g.  on the curvature tensor at the current point $\,\xi\,$, and on the velocity $\,\dot\xi\,$. \  We consider several examples\,:  \parsn
- If the scalar curvature $R(\xi)$ is everywhere non-positive (which is physically relevant, see [L-L]), then the Dirichlet form can be$\,$: \quad  ${\ds -\int_{T^1\M} \Big|\nabla_{\dot\xi}^vF(\xi,\dot\xi) \Big|^2\,R(\xi)\, \mu(d\xi,d\dot\xi)\,}$, \parn
leading to the covariant relativistic diffusion we call \ub{$R$-diffusion}. \parsn
- If the energy $\EE(\xi,\dot\xi)$ is everywhere non-negative (which is physically relevant, see [L-L], and [H-E], where this is called the ``weak energy condition''), then we can choose the Dirichlet form to be$\,$: \quad  ${\ds \int_{T^1\M} \Big|\nabla_{\dot\xi}^vF(\xi,\dot\xi) \Big|^2\,\EE(\xi,\dot\xi)\, \mu(d\xi,d\dot\xi)\,}$, \parn
leading to another covariant relativistic diffusion, we call \ub{energy $\EE$-diffusion}. \parn
Contrary to the basic relativistic diffusion, these new relativistic diffusions reduce to the geodesic flow in every empty (vacuum) region. 
\parsn
Note that  $\, -R(\xi)$ and $\EE(\xi,\dot\xi)$  could be replaced by $\,\f(R(\xi))$ and $\,\psi(\EE(\xi,\dot\xi)$, for more or less arbitrary non-negative increasing functions $\,\f,\psi$. We shall present this class of \ub{covariant $\,\Xi$-relativistic diffusions}, or \ub{$\Xi$-diffusions}, in Section \ref{sec.GCCRD} below. 
\parsn 
- If the sectional curvatures of timelike planes are everywhere non-negative (sectional curvature has proved to be a natural tool in Lorentzian geometry, see for example [H], [H-R]), as this is often the case (at least in usual symmetrical examples), then it is possible to construct a covariant \ub{sectional relativistic diffusion} which undergoes velocity fluctuations that are no longer isotropically Brownian, using the whole curvature tensor (not the Ricci tensor alone). See Section \ref{sec.SRD} below. This sectional relativistic diffusion depends on the curvature tensor in a canonical way. Its diffusion symbol vanishes in flat regions (i.e. regions where the whole curvature tensor vanishes), but does not vanish, in general, in empty regions (i.e. regions where the Ricci tensor vanishes).  \pars 

    Note that all these covariant diffusions are the projections on $T^1\M$ of diffusions on the frame bundle $G(\M)$. Actually, they are constructed directly on $G(\M)$, as in [F-LJ] and in the classical construction of Brownian motion on Riemannian manifolds, see [El], [M], [I-W], [Em], [Hs], [A-C-T].   
These constructions are performed in Sections \ref{sec.GCCRD} and \ref{sec.SRD} below. \par 
   Note also that, while in the flat case the Dudley diffusion [Du] is the unique covariant diffusion, the above examples show that this is not at all the same for curved spaces. \pars  
   
   In Sections \ref{sec.CWP} and \ref{sec.exRW}, we study in more detail the case of warped products, and specify further some particularly symmetrical examples, namely Robertson-Walker manifolds, which are warped products with energy-momentum tensor of perfect fluid type. \par
    We investigate more closely Einstein-de Sitter-like manifolds (Robertson-Walker manifolds for which the expansion rate is $\,\alpha(t) = t^c\,$ for some positive $\,c$), reviewing in this simple class of examples, the relativistic diffusions we introduced, which appear to be distinct. We perform in this setting an asymptotic study of the $R$-diffusion, and of the minimal sub-diffusion $(t_s,\dot t_s)$ relative to the energy $\EE$-diffusion. \par 
   
\section{Canonical vector fields and curvature}  \label{sec.IntCcurvRelDiff} \label{s.CanVFC} 
   We present in this section the main notations and recall a few known facts (see [K-N]). 
   
\subsection{Isomorphism between $\bigwedge^2\bR^{1,d}$ and $\,so(1,d)$}  \label{sec.IsomL^2Rso1d} \indf 
   On the Minkowski space-time $\R^{1,d}$, we denote by $\,\eta = (\!(\eta_{ij})\!)_{0\le i,j\le d}\,$ the Minkowski tensor $\begin{pmatrix} 1& 0 & \ldots & 0\cr 0& -1 & \ldots & 0\cr \vdots &  & \ddots &  \vdots \cr  0& \ldots & 0 & -1 \end{pmatrix}\,$. We also denote by $(\!(\eta^{ij})\!)_{0\le i,j\le d}\,$ the inverse tensor, so that $\,\eta_{ij}\eta^{jk} = \delta _i^k\,$  (or equivalently\,: $\,\eta_{ij} = \eta^{ij} := 1_{\{i=j=0\}} - 1_{\{1\le i=j\le d\}}$), and by $\langle\cdot,\cdot\rangle_\eta$ the corresponding Minkowski pseudo-metric. 
   For $\,u,v,w\in \R^{1,d}$, we set
\begin{equation}  \label{f.isomR1dso} 
u\wedge v\, (w)\, := \, \langle u,w\rangle_\eta\, v - \langle v,w\rangle_\eta\, u\,. 
\end{equation}
In other terms, this is the interior product of $\,u\wedge v\,$ by the dual of $\,w\,$ with respect to $\,\eta\,$. \parn 
This defines an endomorphism of $\,\R^{1,d}$ which belongs to $\,so(1,d)$, since for any $\,w,w'\in \R^{1,d}$ we have clearly \  $\langle u\wedge v (w),w'\rangle_\eta + \langle w,u\wedge v(w')\rangle_\eta = 0\,$. It vanishes only if $\,u\,$ and $\,v\,$ are collinear, hence \iff $\,u\wedge v=0\,$. We have thus an isomorphism between $\bigwedge^2\R^{1,d}$ and $\,so(1,d)$. 
\brem  \label{lem.crochet} \  The Lie bracket of $so(1,d)$ can be expressed, for any $\,a,b,u,v\in \R^{1,d}$, by$\,$: 
\begin{equation*}
[a\wedge b,u\wedge v] = \langle a,u\rangle_\eta\; b\wedge v + \langle b,v\rangle_\eta\; a\wedge u - \langle a,v\rangle_\eta\; b\wedge u - \langle b,u\rangle_\eta\; a\wedge v\, . 
\end{equation*}
\erem 

   The  Minkowski pseudo-metric $\langle\cdot,\cdot\rangle_\eta$ extends to $\bigwedge^2\R^{1,d}$, by setting$\,$:  
\begin{equation}   \label{f.etaL2} 
\langle u\wedge v, a\wedge b\rangle_\eta := \langle u, a\rangle_\eta\,\langle v, b\rangle_\eta - \langle u, b\rangle_\eta\,\langle v, a\rangle_\eta = \5 \big( \langle u\wedge v(a), b\rangle_\eta - \langle u\wedge v(b), a\rangle_\eta \big) , 
\end{equation}
so that, if $(e_0,\ldots,e_d)$ is a Lorentz (i.e. pseudo-orthonormal) basis of $(\R^{1,d}, \eta)$, \  then $(e_i\wedge e_j\,|\, 0\le i<j\le d)$ is an orthogonal basis of $(\bigwedge^2\R^{1,d}, \eta)$, such that $\langle e_i\wedge e_j\,, e_i\wedge e_j\rangle_\eta = \eta_{ii}\,\eta_{jj} \if{ = \eta_{i}\,\eta_{j}}\fi\,$. 

\subsection{Frame bundle $\,G(\M)$ over $(\M,g)$}  \label{sec.FramebGM} \indf 
   Let $\,\M\,$ be a time-oriented $C^\infty$ $(1+d)$-dimensional Lorentz manifold, with pseudo-metric $\, g\,$ having signature $(+,-,\ldots,-)$, and let $\,T^1\M$ denote the positive half of the pseudo-unit tangent bundle. Let $\,G(\M)$ be the bundle of direct pseudo-orthonormal frames, with first element in $\,T^1\M\,$, which has its fibers modelled on the special Lorentz group.   Let $\,\pi_1 : u\mapsto\big(\pi(u), e_0(u)\big)$ denote the canonical projection from $\,G(\M)\,$ onto the unit tangent bundle $\,T^1\M\,$, which to each frame $\big(e_0(u),\ldots, e_d(u)\big)$ associates its first vector $\,e_0\,$. \parn   
    We denote by \  $T\M\mathop{\lra}^{\pi_2}\limits \M\,$  the tangent bundle, by $\,\Gamma(TM)$ the set of $C^2$ vector fields on $\M$ (sections of $\pi_2$), 
by $\,G(\M)\mathop{\lra}^{\pi}\limits \M\,$ the frame bundle, by \  $u=\big(\pi(u) ; e_0(u),\ldots, e_d(u)\big)$ the generic element of $\,G(\M)$.
   We extend (\ref{f.isomR1dso}) to a linear action of $\,so(1,d)\equiv \bigwedge^2\R^{1,d}$ on $G(\M)$, by setting$\,$: 
\begin{equation*}   
e_k\wedge e_\ell\, (e_j(u))\, := \, \eta_{jk}\, e_\ell(u) - \eta_{j\ell}\, e_k(u)\,,  \qquad \hbox{for any }\; 0\le j,k,\ell\le d\, , 
\end{equation*}
where $(e_0,\ldots,e_d)$ denotes the canonical basis of $\,\R^{1,d}$. \parn
The action of $\,SO(d)$ on $(e_1,\ldots,e_d)$ induces the identification $\,T^1\M\equiv G(\M)/SO(d)$. \parn 
The right action of $\,so(1,d)$ on $G(\M)$ defines a linear map $\,V\,$ from  $\,so(1,d)$ into vector fields on $G(\M)$ (i.e. sections of the canonical projection of $\,TG(\M)$ on $G(\M)$), such that 
\begin{equation} \label{f.actso1dGMcr} 
[V_{a\wedge b} ,V_{\alpha\wedge \beta}] = V_{[a\wedge b , \alpha\wedge \beta]}\, ,\qquad \hbox{for any }\; a\wedge b\, , \alpha\wedge \beta\in \bigwedge\!\,\!^2\,\R^{1,d} . 
\end{equation}
Vector fields $\,V_{a\wedge b}\,$ are called  \ub{vertical}. \parn 
{\bf Notation} \  To abreviate the notations, we shall consider mostly the canonical vector fields$\,$:
\begin{equation*}
V_{ij} := V_{e_i\wedge e_j}\, , \quad \hbox{ for }\; 0\le i,j\le d\, .  
\end{equation*}
By (\ref{f.actso1dGMcr}) and (\ref{f.isomR1dso}), for $\,0\le i,j,k,\ell\le d\,$ we have$\,$: 
\begin{equation} \label{f.commrelVVH} 
[V_{ij},V_{k\ell}] = \eta_{ik}\,V_{j\ell} + \eta_{j\ell}\, V_{ik} - \eta_{i\ell}\, V_{jk} - \eta_{jk}\, V_{i\ell} \, . 
\end{equation}
We shall often write $\, V_{j}\,$ for $\,V_{0 j}\, $. \pars  

Denote by $\,p\,$ the canonical projection $\,TT\M\buildrel{p}\over{\lra} TM\,$, by $\,\tilde p\,$ the canonical projection $\,TG(\M)\buildrel{\tilde p}\over{\lra} G(\M)$, and consider also the projection $\,TT\M \buildrel{(p,T\pi_2)}\over{\hbox to 13mm {\rightarrowfill}} T\M\oplus T\M\,$, where
\begin{equation*}
T\M\oplus T\M := \big\{ (\xi\, ; v_1,v_2)\,\big|\, \xi\in \M\,,v_i\in T_\xi\M \big\} \equiv \big\{ (w_1,w_2)\,\big|\, w_i\in T\M \, ,\, \pi_2(w_1) = \pi_2(w_2)\big\}  
\end{equation*}
is the so-called Whitney sum. 

   A \ub{connection} $\,\sigma\,$ can be defined as a bilinear section $\,T\M\oplus T\M \buildrel{\sigma}\over{\lra} TT\M\,$ of $\,(p,T\pi_2)$, the bilinearity being that of $(v_1,v_2)\mapsto \sigma(\xi\, ; v_1,v_2)$, above any given base point $\,\xi\in \M$.  \parn 
Given such a connection $\,\sigma\,$ and a $C^1$ curve $\,\gamma\,$, the \ub{parallel transport} $\,\big/\!\!\big/^{\gamma}_t$ along $\,\gamma\,$ of any $\,v_0\in T_{\gamma_0}\M\,$ is $\,v_t = \big/\!\!\big/^{\gamma}_tv_0\in T_{\gamma_t}\M\,$ defined by the ordinary differential equation$\,$: \   $\frac{d}{dt}(\gamma_t\,; v_t) = \sigma(\gamma_t\,;\,v_t\,, \dot\gamma_t)\,$. \  
Then the \ub{covariant derivative} $\,\nabla_{\dot\gamma_0}X(\gamma_0)\in T_{\gamma_0}\M\,$ of a $C^1$ vector field $X$ is defined as the derivative at 0\, of\, $\,t\mapsto \big(\big/\!\!\big/^{\gamma}_t\big)\1 X(\gamma_t)\,$. \  
\par 
   
     A connection $\,\sigma\,$ is said to be \ub{metric} if the associated parallel transport preserves the pseudo-metric, and then acts on $G(\M)$ as well as on $\,T\M\,$. \par 
 \par
     A metric connection $\,\sigma\,$ defines the \ub{horizontal vector fields} $\,H_k\,$ on $G(\M)$, for $0\le k\le d\,$, given  for any $\,F\in C^1(G(\M))$ and $\,u\in G(\M)$, by$\,$: 
\begin{equation} \label{f.defHe} 
H_kF(u) \; \hbox{ is the  derivative at 0\, of } \;\, t\mapsto F\big(\big/\!\!\big/^{\gamma}_t u\big) , 
\end{equation}
the $C^1$ curve $\gamma\,$ being such that $\,\gamma_0= \pi(u) , \dot\gamma_0= e_k(u)\,$. Note that $\,T\pi(H_k) = e_k\,$. \pars  
The canonical vectors $\,V_{ij}, H_k\,$ span $TG(\M)$ \big(the horizontal (resp. vertical) sub-bundle of $TG(\M)$ being spanned by $\,H_k$'s (resp. $V_{ij}$'s)\big). \    Note that $\,H_0\,$ generates the geodesic flow, that $\,V_1,\ldots,V_d\,$ generate the \ub{boosts}, and that the $\,V_{ij}\, (1\le i,j\le d)$ generate rotations. \parsn 
This allows to define the intrinsic \ub{torsion} tensor $(\!(\TT_{ij}^k)\!)$ and \ub{curvature} tensor $(\!(\RR_{ij}\,\!^{k\ell})\!)$ (with $\,0\le i,j,k,\ell\le d$) of the metric connection $\,\sigma$, by the assignment$\,$: 
\begin{equation} \label{f.defRRTT} 
[H_i,H_j] \, = \sum_{k=0}^d\TT_{ij}^k\,H_{k} + \sum_{0\le k<\ell\le d} \RR_{ij}\,\!^{k\ell}\,V_{k\ell}
\end{equation} we can denote more simply by $\,\TT_{ij}^k\,H_{k} + \5\,\RR_{ij}\,\!^{k\ell}\,V_{k\ell}\, $.  
\parn  
For any metric connection we have$\,$: 
\begin{equation} \label{f.ComrelVH} 
 [V_{ij},H_{k}] = \eta_{ik}\,H_{j} - \eta_{jk}\,H_{i}\, , \qquad \hbox{for }\; 0\le i,j,k\le d\,. 
\end{equation}

   There exists a unique metric connection with vanishing torsion, called the \parn 
\ub{Levi-Civita connection}.  We shall henceforth consider this one. \parsn 
\if{ 
Given two vector fields $X,Y$ on $\M$ and $\,\xi\in\M$, the \ub{covariant derivative} $\,\nabla_XY(\xi)\in T_\xi\M\,$ is 
defined by setting \big(using the base $(e^0(u),\ldots,e^d(u))$ dual to $(e_0(u),\ldots,e_d(u))$\big): 
\begin{equation}  \label{f.nabla} 
\nabla_{X} Y(\xi) :=\, \sum_{0\le k,\ell\le d} \big\langle X(\xi) , e^k(u)\big\rangle\, H_k(u)\big(\langle\, Y\circ\pi , e^\ell\,\rangle\big)\, e_\ell(u) \, , \quad \hbox{ for any }\; u\in\pi\1(\xi)\, . 
\end{equation}
Note that this makes indeed sense, since $\langle\, Y\circ\pi , e^\ell\,\rangle$ is a function on $G(\M)$, on which the tangent vector $\,H_k(u)$ does act, and that $\,\nabla_XY\,$ is bilinear with respect to $(X,Y)$. \pars 
}\fi  
The \ub{curvature operator} $\,\RR_\xi\,$ is defined on $\bigwedge^2T_{\xi}\M\,$ by\,:  
\begin{equation} \label{f.RiemR'} 
\RR_\xi\big(e_i(u)\wedge e_j(u)\big) := \sum_{0\le k<\ell\le d} \RR_{ij}\,\!^{k\ell}\, e_k(u)\wedge e_\ell(u) , \hbox{ for any}\; u\in\pi\1(\xi)\; \hbox{ and }\; 0\le i,j\le d\, . 
\end{equation}
The curvature operator is alternatively given by$\,$:  for any $\,C^1\!$ vector fields $\, X,Y,Z,A\,$, 
\begin{equation} \label{f.exprlocR} 
\left\langle \RR\left(X\wedge Y\right) \hbox{,}\, A\wedge Z\right\rangle_{\!\eta} = \left\langle \left([\nabla\!_{X}\, , \nabla\!_{Y}] - \nabla\!_{[X,Y]}\right) Z\,,\, A\right\rangle_{\! g} .  
\end{equation}
The \ub{Ricci tensor} and \ub{Ricci operator} are defined, for $\,0\le i,k\le d\,$, by$\,$: 
\begin{equation}  \label{f.RicciR} 
R_{i}^{k} := \sum_{j=0}^d \RR_{ij}\,\!^{kj}\,, \quad \hbox{ and }\quad {\rm Ricci}_\xi\big(e_i(u)\big) := \sum_{k=0}^d R_{i}^{k}\,e_k(u)\,  ,   \; \hbox{ for any }\; u\in\pi\1(\xi)\, . 
\end{equation}
The \ub{scalar curvature} is\,: \quad   ${\ds R := \sum_{k=0}^d R_{k}^{k}\,}$. \parsn
The indexes of the curvature tensor $(\!(\RR_{ij}\,\!^{k\ell})\!)$ and of the Ricci tensor $(\!(R_i^k)\!)$ are lowered or raised by means of the Minkowski tensor $(\!(\eta_{ab})\!)$ and its inverse $(\!(\eta^{ab})\!)$. For example, we have\,: \  $\RR^p\,\!_{jqr} = \RR_{ij}\,\!^{k\ell}\,\eta^{ip}\, \eta_{kq}\,\eta_{\ell r}\,$, \  and \  $R_{ij} = R_{i}^{k}\, \eta_{kj}\,$. 

\brem   \label{symmRt}   \  {\rm  The curvature and Ricci operators and tensors are symmetrical$\,$:  \parsn 
\centerline{$\langle \RR(a\wedge b) , v\wedge w\rangle_\eta = \langle a\wedge b , \RR(v\wedge w)\rangle_\eta\,$, \  and \   $\,\langle {\rm  Ricci}(v) , w\rangle_\eta = \langle v , {\rm  Ricci}(w)\rangle_\eta\,$,}  \parsn 
for any $a,b,v,w\in \R^{1,d}$. \   Equivalently, for $\,0\le i,j,k,\ell\le d\,$:   \  
{$R_{ij}\!\,^{k\ell} = R^{k\ell}\!\,_{ ij}$,  and  \   $R_{ij} = R_{ji}\,$.}
}\erem
The \ub{energy-momentum}  tensor $(\!(T^k_j)\!)$ and operator $T_\xi\,$ are defined as$\,$: 
\begin{equation} \label{f.Einstein} 
T^k_j := R^k_j - \5\, R\, \delta^k_j  \quad \hbox{ and }\quad 
T_\xi := {\rm Ricci}_\xi - \5\, R\,  . 
\end{equation} 
Note that \  ${\ds \sum_{j=0}^d T^j_j = - {\ts \frac{d-1}{2}}\, R\,}$. \quad    
The \ub{energy} at any line-element $(\xi,\dot\xi)\in T^1\M$ is  
\begin{equation} \label{f.Ener} 
\EE(\xi,\dot\xi) := \langle T_\xi(\dot\xi) , \dot\xi\rangle _{g(\xi)}  =\,T_{00}(\xi,\dot\xi)\, .   
\end{equation} 
The last equality is easily derived from (\ref{f.RicciR}) and (\ref{f.Einstein}) since writing  $(\xi,\dot\xi)=(\pi(u), e_0(u))$ for any $\,u\in\pi_1\1(\xi,\dot\xi)$  and $\,T_{ij} = T_{i}^{k}\, \eta_{kj}\,$, we have$\,$: 
\begin{equation*}
\langle T_\xi(\dot\xi) , \dot\xi\rangle _{g(\xi)} = g\big( T_\xi(e_0) , e_0\big) = g(T_0^k\, e_k,e_0) = T_0^k\, \eta_{k0} = T_{00} = T_{00}(\xi,\dot\xi)\, . 
\end{equation*}
The \ub{weak energy condition} (see [H-E]) stipulates that \  $\EE(\xi,\dot\xi)\ge 0\,$ on the whole $\,T^1\M$.   \parm 

   We shall need the following general computation rule. 
\blem  \label{lem.calcViRo} \quad   For $\,0\le i,j,k,\ell,p,q\le d\,$, \  we have$\,$: 
\begin{equation*}
V_{qp}\RR_{ij}\!\,^{k\ell} =  \eta_{qi}\, \RR_{pj}\!\,^{k\ell} - \eta_{ip}\, \RR_{qj}\!\,^{k\ell} + \eta_{qj}\, \RR_{ip}\!\,^{k\ell} - \eta_{jp}\, \RR_{iq}\!\,^{k\ell} + \delta_{q}^{k}\, \RR_{ijp}\!\,^{\ell} - \delta_{p}^{k}\, \RR_{ijq}\!\,^{\ell} - \delta_{q}^{\ell}\, \RR_{ijp}\!\,^{k} + \delta_{p}^{\ell} \RR_{ijq}\!\,^{k} . 
\end{equation*}  
\elem 
\ub{Proof} \quad  Using (\ref{f.RiemR'}) and (\ref{f.isomR1dso}), we have indeed$\,$: 
\parn  \vbox{ 
\begin{equation*}
V_{qp}\, \RR_{ij}\!\,^{k\ell} = V_{qp}\, \langle \RR(e_i\wedge e_j) , e_n\wedge e_m\rangle_\eta\, \eta^{kn}\, \eta^{\ell m} 
\end{equation*}  
\begin{equation*}  
=  \bigg[ \Big\langle \RR\Big(\!(\eta_{qi} e_p -\eta_{ip} e_q)\wedge e_j\!\Big) , e_n\wedge e_{m} \Big\rangle_{\!\eta} + \Big\langle \RR\Big(\! e_i\wedge (\eta_{qj} e_p -\eta_{jp} e_q)\!\Big) , e_n\wedge e_{m} \Big\rangle_{\!\eta} \bigg] \eta^{kn} \eta^{\ell m}
\end{equation*}  
\begin{equation*}  
+\, \bigg[ \Big\langle \RR(e_i\wedge e_j) , (\eta_{qn} e_p -\eta_{np} e_q)\wedge e_{m} \Big\rangle_{\!\eta} + \Big\langle \RR(e_i\wedge e_j) , e_n\wedge (\eta_{qm} e_p -\eta_{mp} e_q) \Big\rangle_{\!\eta}\bigg] \eta^{kn} \eta^{\ell m} 
\end{equation*} }\parn 
\begin{equation*}  
=  \eta_{qi}\, \RR_{pj}\!\,^{k\ell} - \eta_{ip}\, \RR_{qj}\!\,^{k\ell} + \eta_{qj}\, \RR_{ip}\!\,^{k\ell} - \eta_{jp}\, \RR_{iq}\!\,^{k\ell} + \delta_{q}^{k}\, \RR_{ijp}\!\,^{\ell} - \delta_{p}^{k}\, \RR_{ijq}\!\,^{\ell} - \delta_{q}^{\ell}\, \RR_{ijp}\!\,^{k} + \delta_{p}^{\ell} \RR_{ijq}\!\,^{k} . \;\;\diamond 
\end{equation*}

\subsection{Expressions in local coordinates}  \label{sec.ExprLocCoord} \indf 
     Consider local coordinates $(\xi^i,e_j^k)$ for $\, u= (\xi,e_0,\ldots, e_d)\in G(\M)$, with $\, {\ds e_j=e_j^k\,\frac{\p}{\p \xi^k}}\,$\raise 1.5pt\hbox{.} \parn    
Then the horizontal and vertical vector fields $\,V_{ij},V_j,H_k\,$, which satify the commutation relations (\ref{f.commrelVVH}),(\ref{f.defRRTT}),(\ref{f.ComrelVH}) of the preceding section \ref{sec.FramebGM},  read as follows.  First, denoting by $\,\Gamma^\ell_{kj}=\Gamma^\ell_{jk}\,$ the \ub{Christoffel coefficients} of the Levi-Civita connexion $\nabla$, we have for $\,0\le i,j\le d\,$:   
\beq \label{f.HV} 
\nabla\!_{\frac{\p}{\p \xi^i}}\, \frac{\p}{\p \xi^j} = \Gamma^k_{ij}(\xi)\,\frac{\p}{\p \xi^k}\, \quad \hbox{ and }  \quad  H_j= e^k_j\,\frac{\p}{\p\xi^k} - e^k_j\,e_i^m\,\Gamma^\ell_{km}(\xi)\, \frac{\p}{\p e^\ell_i}\, \raise 2pt\hbox{.}  \eeq  
This is consistent with (\ref{f.defHe}). Indeed, we have a priori an expression $\,H_j= a^k_j\,\frac{\p}{\p\xi^k} + b^\ell_{ji}\, \frac{\p}{\p e^\ell_i}\,$\raise 2pt\hbox{,} with on one hand $\, a^k_j = \langle T\pi(H_j) , d\xi^k\rangle = \langle e_j , d\xi^k\rangle = e_j^k\,$. On the other hand, for a $C^1$ curve $\gamma\,$ satisfying $\,\gamma_0= \xi\, , \dot\gamma_0= e_j\,$, denoting by $\,e\1$ the matrix inverse to $\,e\equiv (\!(e^k_i)\!)$ we have$\,$: 
\begin{equation*}
e_j^k\,\Gamma^\ell_{km}\,\frac{\p}{\p\xi^\ell} = \nabla_{\dot\gamma_0}\,\frac{\p}{\p\xi^m} = \nabla_{\dot\gamma_0}\big( (e\1)_{m}^i\, e_i\big) = \frac{d_o}{dt}\Big[ \big(\big/\!\!\big/^{\gamma}_t\big)\1 \big( (e\1)_{m}^i\, e_i\big)\big(\big/\!\!\big/^{\gamma}_t u\big)\Big] 
\end{equation*} 
\begin{equation*}
= \frac{d_o}{dt}\Big[ (e\1)_{m}^i \big(\big/\!\!\big/^{\gamma}_t u\big)\Big] e_i
= H_j \big((e\1)_{m}^i\big) \,e_i^\ell \,\frac{\p}{\p\xi^\ell} = - (e\1)_{m}^i\, H_j(e_i^\ell) \,\frac{\p}{\p\xi^\ell} 
\end{equation*}
by (\ref{f.defHe}), so that \  
$b^\ell_{ji} = H_j(e^\ell_i) = - e^k_j\,e_i^m\,\Gamma^\ell_{km}\,$ as wanted. \pars 

   Recall that the Christoffel coefficients of $\nabla$ are computed by$\,$: \par 
\begin{equation*}
\G_{ij}^k = \5 \, g^{k\ell}\Big( \frac{\p g_{\ell j}}{\p\xi^i} + \frac{\p g_{i \ell}}{\p\xi^j} - \frac{\p g_{i j}}{\p\xi^\ell}\Big)  \, ,
\end{equation*}
or equivalently, by the fact that geodesics solve \  $\ddot\xi^k +\Gamma^k_{ij}\dot\xi^i\dot\xi^j = 0\,$. \parsn 
Then for $\,0\le i,j, k\le d\,$: 
\begin{equation*}
V_{ij}\,e_k = \eta_{ik}\, e_j - \eta_{jk}\, e_i = ( \eta_{ik}\, e_j^\ell - \eta_{jk}\, e_i^\ell)\, \frac{\p}{\p \xi^\ell} = (\eta_{iq}\, e_j^m - \eta_{jq}\, e_i^m)\, \frac{\p}{\p e_q^m}\, e_k^\ell\, \frac{\p}{\p \xi^\ell}\, \raise 2pt\hbox{,} 
\end{equation*}
whence for $\,0\le i,j\le d\,$: 
\begin{equation} \label{f.V} 
V_{ij}\, =\, (\eta_{iq}\, e_j^m - \eta_{jq}\, e_i^m)\, \frac{\p}{\p e_q^m}\,  \raise 2pt\hbox{,} 
\end{equation}
that is$\,$: \quad ${\ds V_{ij}= e^k_i\, \frac{\p}{\p e^k_j} - e^k_j\, \frac{\p}{\p e^k_i} }\;$  \  and \quad   ${\ds V_{j}= e^k_0\, \frac{\p}{\p e^k_j} + e^k_j\, \frac{\p}{\p e^k_0} }\;$, \   for $\,1\le i,j\le d\,$. \pars

   The curvature operator expresses in a local chart as$\,$: \  for $\,0\le m,n,p,q\le d\,$, 
\begin{equation}  \label{f.exprlocR'} 
\widetilde \RR_{mnpq} := \Big\langle \RR\left({\ts\frac{\p}{\p \xi^m}}\wedge {\ts\frac{\p}{\p \xi^n}}\right) \raise 1.6pt\hbox{,}\, {\ts\frac{\p}{\p \xi^p}}\wedge {\ts\frac{\p}{\p \xi^q}}\Big\rangle_{\!g} =  g_{mr} \bigg( \Gamma^r_{ps}\,\Gamma^s_{nq} - \Gamma^r_{qs}\,\Gamma^s_{np} + {\frac{\p \Gamma^r_{nq}}{\p \xi^p}} - {\frac{\p \Gamma^r_{np}}{\p \xi^q}}\bigg) . 
\end{equation}
Then, the Ricci operator can be computed similarly, as$\,$: \   for $\,0\le m,p\le d\,$, 
\begin{equation}   \label{f.exprlocRic} 
\tilde R_{mp} := \Big\langle {\rm Ricci}({\ts\frac{\p}{\p \xi^m}})\, \raise 1.6pt\hbox{,}\, {\ts\frac{\p}{\p \xi^p}} \Big\rangle_{\!g}  = \widetilde \RR_{mnpq}\, g^{nq} 
= \Gamma^n_{nq}\Gamma^q_{mp} - \Gamma^n_{pq}\, \Gamma^q_{mn} + {\frac{\p \Gamma^n_{mp}}{\p \xi^n}} - {\frac{\p \Gamma^n_{mn}}{\p \xi^p}}\, \raise 2pt \hbox{.} 
\end{equation}
The scalar curvature and the energy-momentum operator can be computed by\,: 
\begin{equation} \label{f.Einstein'} 
R =\tilde R_{ij}\,g^{ij} \quad \hbox{ and }\quad  \tilde T_{\ell m} = \tilde R_{\ell m} - \5\,R\, g_{\ell m}\, .
\end{equation} 

  To summarize, the Riemann curvature tensor $(\!(\RR_{ij}\,\!^{k\ell})\!)$ is made of the coordinates of the curvature operator $\,\RR\,$ in an orthonormal moving frame, and its indexes are lowered or raised by means of the Minkowski tensor $(\!(\eta_{ab})\!)$, while the curvature tensor $(\!(\tilde R_{mnpq})\!)$ is made of the coordinates of the curvature operator in a local chart, and its indexes are lowered or raised by means of the metric tensor $(\!(g_{ab})\!)$. \parn
To go from one tensor to the other, note that by (\ref{f.exprlocR'}) and (\ref{f.RiemR'}) we have \parn  
$ \RR\big({\ts\frac{\p}{\p \xi^m}}\wedge {\ts\frac{\p}{\p \xi^n}}\big) = \5\,\widetilde \RR_{mn}\,\!^{ab}\, {\ts\frac{\p}{\p \xi^a}}\wedge {\ts\frac{\p}{\p \xi^b}}\,$\raise1.3pt \hbox{,}  \   whence\,: \  
$ e_i^k\,e_j^\ell\, \widetilde \RR_{k\ell}\,\!^{pq} = \RR_{ij}\,\!^{mn}\,e_m^p\,e_n^q\, $, \   or equivalently\,: 
\begin{equation} \label{f.RtildeR} 
\RR_{ijab} =  \widetilde {\RR}_{k\ell rs}\,e_i^k\,e_j^\ell\,e_a^r\,e_b^s\, ,  \quad  \hbox{ or as well\,: } \quad  \widetilde{\RR}^{rspq} = \RR^{abmn}\,e_a^r\,e_b^s\,e_m^p\,e_n^q\, .
\end{equation}

\subsection{Case of a perfect fluid} \label{sec.perffluid} \indf 
   The energy-momentum tensor $\,T$ (of (\ref{f.Einstein}), or equivalently $\tilde T$, recall (\ref{f.Einstein'})) is associated to a \ub{perfect fluid} (see [H-E]) if it has the form :
\begin{equation}  \label{f.fp} 
 \tilde T_{k\ell} = q\, U_k\,U_\ell - p\, g_{k\ell}\, , 
\end{equation}  
for some $\,C^1$ field $\,U\,$ in $\,T^1\M$ (which represents the velocity of the fluid), and some $\,C^1$ functions $\,p,q\,$ on $\,\M\,$.  
By Einstein Equations (\ref{f.Einstein'}),  (\ref{f.fp}) is equivalent to$\,$:  
\begin{equation}  \label{f.fpb}  
\tilde R_{k\ell} = q\, U_k\,U_\ell + \tilde p\, g_{k\ell}\, , \quad \hbox{ with } \quad \tilde p = ({2 p - q})/({d-1}) , 
\end{equation}  
or as well, by (\ref{f.exprlocRic}), to$\,$:  
\begin{equation}  \label{f.fpb'}  
\langle{\rm Ricci}(W) , W\rangle_{\eta} = q\times g(U,W)^2 + \tilde p\times g(W,W)\, , \quad \hbox{ for any } \; W\in T\M\, .   
\end{equation}
The quantity $\,\langle U(\xi_s),\dot\xi_s\rangle$, (which is the hyperbolic cosine of the distance, on the unit hyperboloid at $\,\xi_s\,$ identified with the hyperbolic space between the space-time velocities of the fluid and the path) will be denoted by $\,  \AA_s\,$ or  $\,\AA(\xi_s,\dot\xi_s)$. \parn 
Note that necessarily \  $\AA_s\ge 1\,$.  
By Formulas (\ref{f.Ener}) and (\ref{f.fp}), the energy equals$\,$: 
\begin{equation} \label{f.Enerpf} 
\EE(\xi,\dot\xi) = q(\xi)\, \AA(\xi,\dot\xi)^2 - p(\xi) .  
\end{equation} 

\brem \label{rem.RTpf} \   {\rm     $(i)$ \  The energy of the fluid is simply$\,$:  \  $\tilde T_{k\ell}\,U^k\, U^\ell = q-p\, $. \quad   and the scalar curvature equals \quad  $R = 2\, [{(d+1)\, p - q}]/({d-1}) $. \parn    
$(ii)$ \  By (\ref{f.Enerpf}), the weak energy condition reads here$\,$: \  $q\ge p^+$. 
}\erem

\section{Covariant $\,\Xi$-relativistic diffusions}  \label{sec.GCCRD}  \indf 
   Let $\,\Xi \,$ denote a non-negative smooth function on $G(\M)$, invariant under the right action of $\,SO(d)$ (so that it identifies with a function on $T^1\M$). \parn  
   Our examples will be \quad  $\Xi = -\,\rr^2\, R$ \  and  \quad  $\Xi = \rr^2\,\EE$ \quad   (for a positive constant $\,\rr$). \pars
   
           We call \  \ub{$\Xi$-relativistic diffusion} or \ub{$\Xi$-diffusion} the $\,G(\M)$-valued diffusion process $(\Phi_s)$ associated to $\,\Xi\,$ we will construct in Section \ref{sec.XiGCCRD} below, as well as its $\,T^1\M$-valued projection $\,\pi_1(\Phi_s)$. Let us first recall our previous construction, which corresponds to a constant $\,\Xi\,$. \par 
           
           \subsection{The basic relativistic diffusion} \label{sec.reldiff} \indf 
   The relativistic diffusion process $(\xi_s,\dot\xi_s)$ was defined in [F-LJ], as the projection under $\,\pi_1\,$ of the $\,G(\M)$-valued diffusion $(\Psi_s)$ solving the following Stratonovitch stochastic differential equation (for a given $\,\R^d$-valued Brownian motion $(w^{j}_s)$ and some fixed $\,\rr >0$)\,:
\beq \label{f.reldiff} d\Psi_s = H_0(\Psi_s)\, ds + \rr \sum_{j=1}^dV_{j}(\Psi_s)\circ dw^{j}_s\, . \eeq  
The infinitesimal generator of the  $\,G(\M)$-valued relativistic diffusion $(\Psi_s)$ is 
\begin{equation}  \label{f.basrelop}
\HH\,  :=\, H_0+{\ts\frac{\rr^2}{2}}\sum_{j=1}^d V_j^2\, , 
\end{equation}
and the infinitesimal generator of the relativistic diffusion $(\xi_s,\dot\xi_s) := \pi_1(\Psi_s)$ is the relativistic operator$\,$: 
\beq \label{f.relop}   \HH^1 :=  \LL^0 + {\ts\frac{\rr^2}{2}}\,\Delta^v  =\,  \dot\xi^k \frac{\p}{\p\xi^k} + \Big( {\ts\frac{d\,\rr^2}{2}}\, \dot\xi^k - \dot\xi^i\dot\xi^j\, \G_{ij}^k(\xi) \Big)  \frac{\p}{\p\dot\xi^k} + {\ts\frac{\rr^2}{2}}\, (\dot\xi^k \dot\xi^\ell - g^{k\ell}(\xi)) \,\frac{\p^2}{\p \dot\xi^k \p \dot\xi^\ell} \,\raise 1.5pt\hbox{,} \eeq 
where $\,\LL^0$ denotes the vector field on $T^1\M$ generating the geodesic flow, and $\,\Delta^v\,$ denotes the vertical Laplacian, i.e. the Laplacian on $T^1_\xi\M$ equipped with the hyperbolic metric induced by $\,g(\xi)$. The relativistic diffusion $(\xi_s,\dot\xi_s)$ is parametrized by proper time $\,s\ge 0\,$, possibly till some positive explosion time. \   
\parn 
In  local coordinates $(\xi^i,e_j^k)$, setting $\,\Psi_s = (\xi^i_s,e_j^k(s))$, Equation $(\ref{f.reldiff})$ becomes locally equivalent to the following system of It\^o equations\,: 
\begin{equation*}
d\xi^k_s = \dot\xi_s^k\, ds = e_0^k(s)\, ds \; ; \;\; d\dot\xi_s^k = -\G_{i\ell}^k(\xi_s)\, \dot\xi_s^i\, \dot\xi_s^\ell\, ds + \rr \sum_{i=1}^d e^k_i(s)\, dw^{i}_s + {\ts\frac{d\,\rr^2}{2}}\, \dot\xi_s^k\, ds \; , \quad \mbox{ and} 
\end{equation*} 
\vspace{-2mm}
\begin{equation*}
de^k_j(s)= -\G_{i\ell}^k(\xi_s)\, e^\ell_j(s)\, \dot\xi_s^i\, ds+ \rr\, \dot\xi_s^k\, dw^{j}_s + {\ts\frac{\rr^2}{2}}\, e^k_j(s)\,ds \; ,  \quad \mbox{ for }\, 1\le j\le d\,,\; 0\le k\le d\; . 
\end{equation*}  \pars  

\brem \label{rem.SVi2E} \  {\rm  We have on $\,T^1\M\,$: 
\begin{equation} \label{f.HEner} 
\sum_{j=1}^d V_j^2 \EE = 2(d+1)\, \EE - 2\,\Tr(T) = 2(d+1)\, \EE + (d-1)\, R\, .   
\end{equation} 
}\erem 
Indeed, since for each $\,1\le j\le d\,$ $V_j\,$ exchanges the basis vectors $\,e_0=\dot\xi\,$ and $\,e_j$ (recall (\ref{f.V})) we get$\,$: $V_j^2 \EE = 2\, V_j(\tilde T_{\ell m}\,\dot\xi^\ell\,e_j^m) = 2\, \tilde T_{\ell m}\,(\dot\xi^\ell\,\dot\xi^m + e_j^\ell\,e_j^m)$, \   whence 
\begin{equation*}
\sum_{j=1}^d V_j^2 \EE = 2d\, \EE + 2\,\tilde T_{\ell m}\,(\dot\xi^\ell\,\dot\xi^m-g^{\ell m}) = 2(d+1)\, \EE - 2\,\Tr(\tilde T) \, , 
\end{equation*}
and Formulas (\ref{f.HEner}) follow at once by (\ref{f.Einstein}). \parn 
As an application, a direct computation yields the following evolution of the energy.  
\brem \label{cor.reldiff} \   The random energy process $\,\EE_s = \EE(\xi_s,\dot\xi_s)$ associated to the basic relativistic diffusion $\,\pi_1(\Psi_s) = (\xi_s,\dot\xi_s)$ satisfies the following equation (where $\,\nabla_v:= v^j\nabla_j$)$\,$: 
\begin{equation*} 
d\EE_s = \nabla_{\dot\xi_s}\EE_s\, ds + \rr^2 \left[(d+1) \EE_s + {\ts\frac{d-1}{2}}\, R(\xi_s)\right] ds +  dM^{\EE}_s\, , 
\end{equation*} 
with the quadratic variation of its martingale part $\,dM^{\EE}_s\,$ given by$\,$: 
\begin{equation*}
[d \EE_s , d \EE_s] = [d M^{\EE}_s , d M^{\EE}_s] = 4\rr^2\, [\EE_s^2 - \langle \tilde T\dot\xi_s , \tilde T\dot\xi_s\rangle]\, ds\, .
\end{equation*} 
\erem  
\centerline{\big(We have here in particular \quad $\nabla_{\dot\xi_s}\EE_s = \left[\partial_{\xi^k} \tilde T_{ij}(\xi_s)\, - 2\,\tilde T_{i\ell}(\xi_s) \, \G_{jk}^\ell(\xi_s)\, \right] \dot\xi_s^i\, \dot\xi_s^j\,\dot\xi_s^k $.\big)} \parn 
Note that the energy $\,\EE_s\,$ is not, in general, a Markov process. \par 

\subsection{Construction of the $\,\Xi$-diffusion}  \label{sec.XiGCCRD}  \indf 
   Let us start with the following Stratonovitch stochastic differential equation on $G(\M)$ (for a given $\,\R^d$-valued Brownian motion $(w^{j}_s)$)\,:
\beq \label{f.Xireldiff} d\Phi_s = H_0(\Phi_s)\, ds + {\ts\frac{1}{4}}  \sum_{j=1}^d V_{j}\,\Xi\, (\Phi_s) V_{j}(\Phi_s)\, ds + \sum_{j=1}^d\sqrt{\Xi(\Phi_s)}\, V_{j}(\Phi_s)\circ dw^{j}_s\, . \eeq  
Note that all coefficients in this equation are clearly smooth, except $\,\sqrt{\Xi}\,$ on its vanishing set $\,\Xi\1(0)$. However, $\,\sqrt{\Xi}\,$ is a locally Lipschitz function\,; see ([I-W], Proposition IV.6.2).  Hence,  Equation (\ref{f.Xireldiff}) does define a unique $\,G(\M)$-valued diffusion $(\Phi_s)$. \parn  
   We have the following proposition, defining the \,\ub{$\Xi$-relativistic diffusion} (or\, 
\ub{$\Xi$-diffusion}) $(\Phi_s)$ on $\,G(\M)$, and $(\xi_s,\dot\xi_s)$ on $T^1\M$, possibly till some positive explosion time.  
\bpro  \label{the.gen} \  The Stratonovitch stochastic differential equation (\ref{f.Xireldiff}) has a unique solution $\,(\Phi_s) = (\xi_s\,; \,\dot\xi_s\,,e_1(s),\ldots,e_d(s))$, possibly defined till some positive explosion time. This is a $\,G(\M)$-valued covariant diffusion process, with generator 
\begin{equation} \label{f.Xirelop}
\HH_\Xi\, :=\, H_0+ {\ts \frac{1}{2}}\,\sum_{j=1}^{d}\limits V_j\, \Xi\,V_j\, . 
\end{equation}
Its projection $\,\pi_1(\Phi_s) = (\xi_s,\dot\xi_s)$  defines a covariant diffusion on $T^1\M$, with $SO(d)$-invariant generator 
\begin{equation} \label{f.XirelopT1}
\HH_\Xi^1\, :=\, \LL^0 + {\ts \frac{1}{2}}\, \nabla^v\, \Xi\,\nabla^v\,, 
\end{equation}
$\nabla^v\,$ denoting the gradient on $T^1_\xi\M$ equipped with the hyperbolic metric induced by $\,g(\xi)$. \par  
   Moreover, the adjoint of $\,\HH_\Xi\,$ with respect to the Liouville measure of $G(\M)$ is  \parn  
$\HH_\Xi^* :=\, - H_0+ {\ts \frac{1}{2}}\,\sum_{j=1}^{d}\limits V_j\, \Xi\,V_j\,$. \   In particular, if there is no explosion, then the Liouville measure is invariant. Furthermore, if $\,\Xi\,$ does not depend on $\,\dot\xi\,$, i.e. is a function on $\M$, then the Liouville measure is preserved by the stochastic flow defined by Equation (\ref{f.Xireldiff}). 
\epro 
We specify at once how this looks in a local chart, before giving a proof for both statements. 
\bcor  \label{pro.gen} \    The $\,T^1\M$-valued $\,\Xi$-diffusion $(\xi_s,\dot\xi_s)$ satisfies \   $d\xi_s=\dot\xi_s\, ds\,$, \  and  in any local chart, the following It\^o stochastic differential equations$\,$: \quad for $\; 0\le k\le d\,$, (denoting $\,\Xi_s = \Xi(\xi_s,\dot\xi_s)$) \parsn
\begin{equation} \label{f.dxipoint}
d\dot\xi_s^k = dM^k_s - \G_{ij}^k(\xi_s)\, \dot\xi_s^i\, \dot\xi_s^j\, ds + {\ts \frac{d}{2}}\, \Xi_s\, \dot\xi_s^k\, ds + {\ts \frac{1}{2}}\, [\dot\xi^k_s\,\dot\xi^\ell_s - g^{k\ell}(\xi_s)]\,\frac{\partial \Xi}{\partial \dot\xi^\ell}(\xi_s,\dot\xi_s)\, ds \, , 
\end{equation} 
with the quadratic covariation matrix of the martingale term $(dM_s)$ given by$\,$: 
\begin{equation*}
[d \dot\xi^k_s , d \dot\xi^\ell_s] =  [\dot\xi^k_s\, \dot\xi^\ell_s- g^{k \ell}(\xi_s)]\, \Xi_s\, ds\, , \quad \hbox{for }\; 0\le k,\ell\le d\,.
\end{equation*} 
\ecor 
\ub{Proof} \quad  In local coordinates  $(\xi^i,e_j^k)$,  $\,\Phi = (\xi\,;\, e_0\,,e_1,\ldots,e_d)$, using Section \ref{sec.ExprLocCoord}, Equation (\ref{f.Xireldiff}) reads\,: \  for any $\,0\le k\le d\,$, \quad $d\xi^k_s = \dot\xi_s^k\, ds = e_0^k(s)\, ds\, $,  
\begin{equation*}
 d\dot\xi_s^k = -\G_{i\ell}^k(\xi_s)\, \dot\xi_s^i\, \dot\xi_s^\ell\, ds + {\ts\frac{1}{4}} \sum_{j=1}^d V_j\,\Xi\,(\xi_s, \dot\xi_s)\, e^k_j(s)\, ds + \sum_{j=1}^d \sqrt{\Xi(\xi_s, \dot\xi_s)}\,e^k_j(s)\circ dw^{j}_s\; ;
\end{equation*} 
and for \  $ 1\le j\le d\,,\; 0\le k\le d\,$, 
\begin{equation*}
de^k_j(s)= -\G_{i\ell}^k(\xi_s)\, e^\ell_j(s)\, \dot\xi_s^i\, ds + {\ts\frac{1}{4}}\, V_j\,\Xi\,(\xi_s, \dot\xi_s)\, \dot\xi_s^k\,ds + \sqrt{\Xi(\xi_s, \dot\xi_s)}\, \dot\xi_s^k\circ dw^{j}_s\;  . 
\end{equation*}  \parsn
We compute now the It\^o corrections, which involve the partial derivatives of $\sqrt{\Xi(\xi, \dot\xi)}\,$    
with respect to $\,\dot\xi\,$. We get successively, for $\, 1\le j\le d\,,\; 0\le k\le d\,$: \parn 
\begin{equation*}
\sqrt{\Xi(\xi_s, \dot\xi_s)} \left[d\left(\!\sqrt{\Xi(\xi_s, \dot\xi_s)}\,e^k_j(s)\!\right)\! , dw^{j}_s\right] = \, \Xi(\xi_s, \dot\xi_s) \left[d e^k_j(s) , \, dw^{j}_s\right] + \5\, e^k_j(s) \left[d\,\Xi(\xi_s, \dot\xi_s) , dw^{j}_s\right]
\end{equation*} 
\begin{equation*}
= \, \Xi\,(\xi_s, \dot\xi_s)^{3/2}\, \dot\xi_s^k\, ds + \5\, e^k_j(s)\, \frac{\partial\,\Xi}{\partial\dot\xi^\ell}(\xi_s, \dot\xi_s) \sqrt{\Xi(\xi_s, \dot\xi_s)}\, e^\ell_j(s)\, ds \, , 
\end{equation*} 
hence, 
\begin{equation*}
\left[d\left(\!\sqrt{\Xi(\xi_s, \dot\xi_s)}\,e^k_j(s)\!\right)\! , dw^{j}_s\right] = \, \Xi\,(\xi_s, \dot\xi_s)\, \dot\xi_s^k\, ds + \5\, V_j\,\Xi\,(\xi_s, \dot\xi_s)\,e^k_j(s)\, ds \, ; 
\end{equation*} 

\begin{equation*}
\sqrt{\Xi(\xi_s, \dot\xi_s)} \left[d\left(\!\sqrt{\Xi(\xi_s, \dot\xi_s)}\,\dot\xi_s^k\!\right)\! , dw^{j}_s\right] = \, \Xi(\xi_s, \dot\xi_s) \left[d \dot\xi_s^k\, , \, dw^{j}_s\right] + \5\, \dot\xi_s^k \left[d\,\Xi(\xi_s, \dot\xi_s) , dw^{j}_s\right]
\end{equation*} 
\begin{equation*}
= \, \Xi\,(\xi_s, \dot\xi_s)^{3/2}\, e^k_j(s)\, ds + \5\, \dot\xi_s^k\, \frac{\partial\,\Xi}{\partial\dot\xi^\ell}(\xi_s, \dot\xi_s) \sqrt{\Xi(\xi_s, \dot\xi_s)}\, e^\ell_j(s)\, ds \, , 
\end{equation*} 
hence, 
\begin{equation*}
\left[d\left(\!\sqrt{\Xi(\xi_s, \dot\xi_s)}\,\dot\xi_s^k\!\right)\! , dw^{j}_s\right] =\, \Xi\,(\xi_s, \dot\xi_s)\, e^k_j(s)\, ds + \5\, V_j\,\Xi\,(\xi_s, \dot\xi_s)\,\dot\xi_s^k\, ds \, . 
\end{equation*} 
Note that the simplification by $\sqrt{\Xi(\xi_s, \dot\xi_s)} \,$ is allowed, since on $\,\Xi\1(0)$ both sides of the simplified formula vanish identically. \parn 
Hence, in local coordinates and in It\^o form, Equation (\ref{f.Xireldiff}) reads\,: \  for any $\,0\le k\le d\,$, \quad $d\xi^k_s = \dot\xi_s^k\, ds = e_0^k(s)\, ds\, $,  \   and setting \  $dM^k_s := \sum_{j=1}^d\limits \sqrt{\Xi(\xi_s, \dot\xi_s)}\,e^k_j(s)\, dw^{j}_s\,$, 
\begin{equation*}
 d\dot\xi_s^k = dM^k_s  -\G_{i\ell}^k(\xi_s)\, \dot\xi_s^i\, \dot\xi_s^\ell\, ds + {\ts\frac{d}{2}}\, \Xi\,(\xi_s, \dot\xi_s)\, \dot\xi_s^k\, ds + {\ts \frac{1}{2}}\, [\dot\xi^k_s\,\dot\xi^\ell_s - g^{k\ell}(\xi_s)]\,\frac{\partial \Xi}{\partial \dot\xi^\ell}(\xi_s,\dot\xi_s)\, ds \; ;
\end{equation*} 
\begin{equation*}
de^k_j(s)= \sqrt{\Xi(\xi_s, \dot\xi_s)}\, \dot\xi_s^k\, dw^{j}_s -\G_{i\ell}^k(\xi_s)\, e^\ell_j(s)\, \dot\xi_s^i\, ds + {\ts\frac{1}{2}}\, \Xi\,(\xi_s, \dot\xi_s)\, e^k_j(s)\,ds + {\ts \frac{1}{2}}\, V_j\,\Xi\,(\xi_s, \dot\xi_s)\, \dot\xi_s^k\, ds\,  . 
\end{equation*}  \parn
Note that we used the formula \  $\sum_{j=1}^d\limits e_j^k(s) e_j^\ell(s) = \dot\xi^k_s\,\dot\xi^\ell_s - g^{k\ell}(\xi_s)\,$, which expresses that $\,\Phi\in G(\M)$. \  Using again this formula,  we get the quadratic covariation matrix of the martingale term $(dM_s)$, displayed in the above proposition, which shows that $(\pi_1(\Phi_s))$ is indeed a diffusion, and proves Corollary \ref{pro.gen}. \par 
   On the other hand, comparing the above equations with Equations (\ref{f.reldiff}) and (\ref{f.basrelop}), which correspond to $\,\Xi\equiv \rr^2\,$, we get  precisely the wanted form (\ref{f.Xirelop}) for the generator of $(\Phi_s)$. \parn 
Then, since  \   $\nabla^v\, \Xi\,\nabla^v =\, \Xi\, \Delta^v + (\nabla^v\, \Xi)\,\nabla^v\,$, comparing the above equation (\ref{f.dxipoint}) for $\,\dot\xi^k_s\,$ with  Equation (\ref{f.relop}) (for which $\,\Xi\equiv \rr^2$), we see that establishing Formula (\ref{f.XirelopT1}) giving the projected generator $\,\HH_\Xi^1\,$ reduces now to proving that 
\   ${\ds (\nabla^v\, \Xi)\,\nabla^v \equiv [\dot\xi^k_s\,\dot\xi^\ell_s - g^{k\ell}]\,\frac{\partial \Xi}{\partial \dot\xi^\ell}\, \frac{\partial }{\partial \dot\xi^k}\,}$\raise 1.8pt\hbox{.} \parn 
Now this becomes clear, noting that \   ${\ds\nabla^v_j = e_j^k \frac{\partial }{\partial \dot\xi^k}\,}$\raise 1.5pt\hbox{,} \   for each $\,1\le j\le d\,$. \par

   Finally, the assertions relative to the Liouville measure are direct consequences of the fact that the vectors $\,H_0\,$ and $\,V_j\,$ are antisymmetric with respect to the Liouville measure of $G(\M)$. (The invariance with respect to $\,H_0\,$, i.e. to the geodesic flow, is proved in the same way as in the Riemannian case, and the invariance with respect to $\,V_j\,$ is straightforward.)  $\;\diamond$  

\brem \label{rem.Ricp} \  {\rm  $(i)$ \    The vertical terms could be seen as an effect of the matter or the radiation present in the space-time $\M$. The $\,\Xi$-diffusion $(\Phi_s)$ reduces to the geodesic flow in the regions of the space where $\,\Xi\,$ vanishes, which happens in particular for empty space-times $\M$ in the cases \  $\Xi = - \rr^2\, R(\xi)$,  \  or \  $\Xi = \rr^2\, \EE(\xi,\dot\xi)$, \   or also \   $\Xi = -\rr^2 \, R(\xi)\,e^{\kappa\,\EE(\xi,\dot\xi)/R(\xi)}\,$ (for positive constant $\,\kappa$) for example. \par   
   $(ii)$ \  As well as for the basic relativistic diffusion, the law of the $\,\Xi$-relativistic diffusion is covariant with any isometry of $(\M,g)$.  \   The basic relativistic diffusion corresponds to $\,\Xi\equiv \rr^2>0\,$, and the geodesic flow to $\,\Xi\equiv 0\,$.  \par
   $(iii)$ \  In [B] is considered a general model for relativistic diffusions, which may be covariant or not. Up to enlarge it slightly, by allowing the ``rest frame'' (denoted by ${z}$ in [B]) 
to have space vectors of non-unit norm, this model includes the generic $\Xi$-diffusion
\big(compare the above equation (\ref{f.Xireldiff}) to (2.5),(3.3) in [B]\big). 
}\erem 
\vspace{-4mm}  

\subsection{The $\,R$-diffusion}  \label{sec.RLRD} \indf 
       We assume here that the scalar curvature $\,R=R(\xi)$ is everywhere non-positive on $\,\M$, which is physically relevant (see [L-L]), and consider the particular case of Section \ref{sec.XiGCCRD} corresponding to  \  $\Xi = -\rr^2\, R(\xi)$, with a constant positive parameter $\,\rr\,$. \parn  
In this case, as its central term clearly vanishes, Equation (\ref{f.Xireldiff}) takes on the simple form\,:
\begin{equation*}
d\Phi_s = H_0(\Phi_s)\, ds + \rr \sum_{j=1}^d\sqrt{- R(\Phi_s)}\, V_{j}(\Phi_s)\circ dw^{j}_s\, .
\end{equation*}

\vspace{-4mm}  
        
\subsection{The $\,\EE$-diffusion}  \label{sec.RicciRelD} \indf 
       We assume that the Weak Energy Condition (recall Section \ref{sec.FramebGM}) holds (everywhere on $\,T^1\M$), which is physically relevant (see [L-L], [H-E]), and consider the particular case of Section \ref{sec.XiGCCRD} corresponding to  \  $\Xi = \rr^2\, \EE = \rr^2\, \EE(\xi,\dot\xi) = \rr^2\,T_{00}\,$. \parn 
We call \ub{energy relativistic diffusion} or \,\ub{$\EE$-diffusion} the $\,G(\M)$-valued diffusion process $(\Phi_s)$ we get in this way, as well as its $\,T^1\M$-valued projection $\,\pi_1(\Phi_s)$. \parsn 
The following is easily derived from Lemma \ref{lem.calcViRo} and Formula (\ref{f.RicciR}). As a consequence, the central drift term in Equation (\ref{f.Xireldiff}) is a function of the Ricci tensor alone when $\,\Xi\,$ is.  
\blem \label{lem.calcViR1} \  We have  \   $ V_{j}R_{i}^{k} =  \delta_{0i}\, R_{j}^{k} - \eta_{ij}\, R_{0}^{k} + \delta_{0}^{k}\, R_{ij} - \delta_{j}^{k}\, R_{0i}\,$, for $\,0\le i,k\le d\,$ and $\,1\le j\le d\,$. \   In particular, \   $ V_{j}R = 0\,$, \  and \  $ V_{j}\EE = V_jT_{00} =  V_jR_{00} = 2R_{0j} \,$. \elem 

    By Lemma \ref{lem.calcViR1},  the drift term of Equation (\ref{f.dxipoint}) which involves  the derivatives of $\,\Xi\,$ equals here\,: \  ${\ds \rr^2 \sum_{j=1}^d R_{0j}(\xi_s)\, e_j^k(s)\, ds\,}$.  \quad   As we have $\,R_{0j} = \tilde R_{mn}\, e_0^m\, e_j^n\,$ by (\ref{f.exprlocRic}), \   we get the alternative expression\,: \quad    ${\ds \rr^2 \sum_{j=1}^d R_{0j}(\xi_s)\, e_j^k(s)\, ds = \rr^2\, \tilde R_{mn}(\xi_s)\, \dot\xi^m \,[\dot\xi^k_s\,\dot\xi^n_s - g^{kn}(\xi_s)]\, ds \,}$. \parn 
Another expression is got by using Einstein equation (\ref{f.Einstein'}), or equivalently, by computing directly from (\ref{f.dxipoint}) and (\ref{f.Ener})\,: \  ${\ds [\dot\xi^k \dot\xi^\ell - g^{k\ell}]\,\frac{\partial \EE}{\partial \dot\xi^\ell} = 2 [\dot\xi^k \dot\xi^\ell - g^{k\ell}]\,\tilde T_{\ell m} \dot\xi^m = 2 [\EE\,\dot\xi  - \tilde T\dot\xi]^k }$,
where the notation $\,(\tilde T\dot\xi)^k \equiv \tilde T^k_m\dot\xi^m\,$ has the meaning of a matrix product. \parsn   
Hence, Formula (\ref{f.dxipoint}) of Proposition \ref{pro.gen} expressing $\,d\dot\xi_s\,$ reads here\,: 
\begin{equation} \label{f.dxipEdiff}  
d\dot\xi_s^k = dM^k_s - \G_{ij}^k(\xi_s)\, \dot\xi_s^i\, \dot\xi_s^j\, ds + {\ts \frac{\rr^2d}{2}}\, \EE_s\, \dot\xi_s^k\, ds + {\rr^2}\, \tilde R_{mn}(\xi_s)\, \dot\xi^m \,[\dot\xi^k_s\,\dot\xi^n_s - g^{kn}(\xi_s)]\, ds \, , 
\end{equation} 
or equivalently\,:
\begin{equation} \label{f.dxipEdiff'}  
d\dot\xi_s = dM_s - \G_{ij}^{\bf \cdot}(\xi_s)\, \dot\xi_s^i\, \dot\xi_s^j\, ds + \rr^2( {\ts \frac{d}{2}} + 1)\, \EE_s\, \dot\xi_s\, ds - {\rr^2}\, \tilde T\dot\xi_s\, ds \, . 
\end{equation} 
   
   We can then compute the equation satisfied by the random energy $\,\EE_s\,$. In particular, the drift term discussed above for $\,d\dot\xi_s^k\,$ contributes now for\,: \parn   
\centerline{${\ds 2\rr^2\, \tilde T_{km}\, \dot\xi^m [\EE\,\dot\xi  - \tilde T\dot\xi]^k =  2\rr^2 ( \EE\, \tilde T_{km}\,  \dot\xi^m \dot\xi^k -  \tilde T_{km}\,  \dot\xi^m [ \tilde T\dot\xi]^k) =  2\rr^2 ( \EE^2 - [\tilde T\dot\xi]_k\,  [\tilde T\dot\xi]^k)  }$.} \parsn  
   This leads to the following, to be compared with Corollary \ref{cor.reldiff}. 
\brem \label{cor.Rreldiff} \   The random energy $\,\EE_s := \EE(\xi_s,\dot\xi_s)$ associated to the $\EE$-diffusion $(\Phi_s)$ satisfies the following equation (where $\,\nabla_v:= v^j\nabla_j$)$\,$: 
\begin{equation*} 
d\EE_s = \nabla_{\dot\xi_s}\EE(\xi_s,\dot\xi_s)\, ds + (d+2)\,\rr^2\, \EE_s^2\, ds -  2\rr^2\, g(\tilde T\dot\xi_s , \tilde T\dot\xi_s)\,ds  + 2\rr\, dM^{\EE}_s\, ,   
\end{equation*} 
with the quadratic variation of its martingale part $\,dM^{\EE}_s\,$ given by$\,$: 
\begin{equation*}
[d \EE_s , d \EE_s] = 4\rr^2\, [d M^{\EE}_s , d M^{\EE}_s] = 4\rr^2\, [\EE_s^2 - g(\tilde T\dot\xi_s , \tilde T\dot\xi_s)]\, \EE_s\, ds\, .
\end{equation*} 
\erem 
Note that the diffusion coefficient $[\EE_s^2 - g(\tilde T\dot\xi_s , \tilde T\dot\xi_s)]$ appearing in the quadratic variation $[d \EE_s , d \EE_s]$ is necessarily non-negative, which can of course be checked directly.

\brem  \label{rem.EnDiff=BRD} \  {Case of Einstein Lorentz manifolds.} \par {\rm 
  The Lorentz manifold $\M$ is said to be Einstein if its Ricci tensor is proportional to its metric tensor. Bianchi's contracted identities (see [H-E] or [K-N]), which entail the conservation equations $\,\nabla_k\tilde T^{jk}=0\,$, force the proportionality coefficient $\,\tilde p\,$ to be constant on $\M$. Hence$\,$: 
\begin{equation*}  
\tilde R_{\ell m}(\xi) = \tilde p\; g_{\ell m}(\xi) \, , \quad \hbox{ for any } \, \xi\, \hbox{ in }\, \M \, \hbox{ and }\, 0\le \ell,m\le d\,. 
\end{equation*}  
Then the scalar curvature is $\,R(\xi) = (d+1) \tilde p\,$, and by Einstein Equations (\ref{f.Einstein'}) we have$\,$: 
\begin{equation*}  
 \tilde T_{\ell m}(\xi) = (\Lambda - {\ts\frac{d-1}{2}}\, \tilde p)\, g_{\ell m}(\xi)\, =:\, -p\, g_{\ell m}(\xi) \,  . 
 \end{equation*} 
Hence Equation (\ref{f.fp}) holds, \   with $\,q = 0\,$: we are in a limiting case of perfect fluid.  \   Moreover, $\,R\,$ and $\,\EE\,$ 
are constant, so that on an Einstein Lorentz manifold, the $\,R$-diffusion and  the $\,\EE$-diffusion coincide with the basic relativistic diffusion (of Section \ref{sec.reldiff}). 
}\erem 

\vspace{-4mm}  

\section{Warped (or skew) products} \label{sec.CWP} \indf 
    Let us consider here a Lorentz manifold $(\M,g)$ having the warped product form$\,$: $\,\M= I\times M\,$, where $\,I\,$ is an open interval of $\,\R_+\,$ and $(M,h)$ is a Riemannian manifold, is endowed with the Lorentzian pseudo-norm $\,g\,$ given by$\,$: 
\begin{equation}  \label{f.CWP} 
d s^2 := d t^2 - \alpha(t)^2\, |dx|_h^2 \, , 
\end{equation}  
or equivalently 
\begin{equation}  \label{f.CWP'} 
g_{0k} := \delta_{0k} \quad\hbox{ and }\quad g_{ij} := -\alpha(t)^2\, h_{ij}(x)\,  , \;\hbox{ for }\; 0\le k\le d\, ,\; 1\le i,j\le d\,.
\end{equation}  
Here $\,\xi\equiv (t,x)\in I\times M\,$ denotes the generic point of $\M$, and the expansion factor $\,\alpha\,$ is a positive $C^2$ function on $\,I\,$. The so-called Hubble function is$\,$: 
\begin{equation}  \label{f.CWPa} 
H(t):=\alpha'(t)/\alpha(t) \, . 
\end{equation}  
This structure is considered in [B-E], which contains most of the following proposition. 
\bPro  \label{lem.wpcurvt} \  Consider a Lorentz manifold $(\M,g)$ having the warped product form. \parn 
$(i)$ \  Its curvature operator $\,\RR\,$ is  given by$\,$: \   for $\,u,v,w,a\in C^2(I)$ and $\,X,Y,Z,A\in \Gamma(TM)$, 
\begin{equation*}
\Big\langle \RR\left((u\partial_t + X)\wedge (v\partial_t + Y)\right) \hbox{,}\, (a\partial_t + A)\wedge (w\partial_t + Z)\Big\rangle_{\!g} =\, - \alpha^2 \left\langle \KK\left(X\wedge Y\right),\, A\wedge Z\right\rangle   \qquad 
\end{equation*}
\begin{equation} \label{f.curvRwp} 
\qquad + \alpha \alpha''\, h(u\,Y- v\,X , a\, Z-w\,A) + (\alpha \alpha')^2\, [ h(X,Z)\, h(Y,A) - h(X,A)\, h(Y,Z)]\, ,
\end{equation} 
$\,\KK\,$ denoting the curvature operator of $(M,h)$. \parn
$(ii)$ \   Denoting by $\,{\rm Ric}\,$ the Ricci operator of $(M,h)$ and by $\langle \cdot\hbox{,}\cdot\rangle$ the standard canonical inner product of $\R^d$, we have\,: 
\begin{equation} \label{f.RicRwp} 
\big\langle {\rm Ricci}(v\partial_t + Y) \hbox{,}\, w\partial_t + Z\big\rangle_{\!\eta} = \langle {\rm Ric}(Y) , Z\rangle + [(d-1)\, |\alpha'|^2 + \alpha \alpha'']\, h(Y,Z) - d\, \frac{\alpha''}{\alpha}\, v w . 
\end{equation}  
$(iii)$ \  If $\,\nabla^{g}$ and $\,\nabla^{h}$ denote the Levi-Civita connections of $(\M,g)$ and $(M,h)$ respectively$\,$: 
\begin{equation} \label{f.LCcwp} 
\nabla^{g}_{(u\partial_t + X)} (v\partial_t + Y) = \left[ uv' + \alpha(t)\alpha'(t)\, h(X,Y)\right] \partial_t + \nabla^h_XY + H(t)\, (u\,Y+v\,X)\, . 
\end{equation} 
\ePro 
\bRem \label{rem.wpcurvt}  \   In any warped local chart $(t, x^j)\,$: for $\, 1\le m,n,p,q\le d\,$ ($K_{mnpq}\,$ denoting the curvature tensor of $(M,h)$, and $K_{mp}$ its Ricci tensor) and $\,0\le k\le d\,$,  \  we have 
\begin{equation} \label{f.LCcwploc} 
\Gamma_{00}^{\bf\cdot}(g) =  0\; ; \quad \Gamma_{0n}^p(g) =  H\, \delta_n^p\; ; \quad 
\Gamma_{mn}^0(g) = \alpha \alpha'\, h_{mn}\; ; \quad  \Gamma_{mn}^p(g) =  \Gamma_{mn}^p(h)\, , 
\end{equation}
\begin{equation} \label{f.curvRwploc} 
\widetilde \RR_{0nkq} =  \delta_{0k}\, \alpha(t) \alpha''(t)\,h_{nq}\; ; \quad 
\widetilde \RR_{mnpq} =  (\alpha \alpha')^2(t)\, [h_{mq}\, h_{np} - h_{mp}\, h_{nq}] - \alpha^2(t) \widetilde K_{mnpq}\, ,  
\end{equation}
\begin{equation} \label{f.curvRicwploc} 
\tilde R_{0k} =  -\, \delta_{0k}\, d\times {\alpha''}(t)/{\alpha}(t) \; ; \quad 
\tilde R_{mp} =  \tilde K_{mp} + [(d-1)\, |\alpha'|^2 + \alpha \alpha'']\, h_{mp} \, .  
\end{equation}
\eRem
Let us outline the proof for the convenience of the reader. \parn 
\ub{Proof} \quad  We get (\ref{f.LCcwp}), and then (\ref{f.LCcwploc}), by using Koszul formula$\,$: 
\begin{equation*}
2\, h(\nabla^h_XY,Z) = Xh(Y,Z) + Yh(X,Z) -Zh(X,Y) + h([X,Y],Z) - h([X,Z],Y) + h([Z,Y],X) , 
\end{equation*}
for both $\,\nabla^{h}$ and $\,\nabla^{g}$. \   Hence, 
\begin{equation*}
\nabla^g_{(u\partial_t + X)} \nabla^g_{(v\partial_t + Y)} (w\partial_t + Z) = \nabla^g_{(u\partial_t + X)} \Big( [vw'+\alpha\alpha'\, h(Y,Z)]\partial_t + \nabla^h_YZ + H\, [vZ+wY] \Big)
\end{equation*}
\vspace{-7mm} 
\parn\vbox{ 
\begin{equation*}
= \Big( u [v'w'+vw''+(\alpha\alpha')' h(Y,Z)] + \alpha\alpha' \Big[ X h(Y,Z) + h(X,\nabla^h_YZ) + H [v h(X,Z)+w h(X,Y)]\Big] \Big)\partial_t 
\end{equation*}
\begin{equation*}
+ \nabla^h_X\nabla^h_YZ + H\, [v\,\nabla^h_XZ+w\,\nabla^h_XY] + uH'\, [v\,Z+w\,Y] + uH\, [v'\,Z+w'\,Y]\, , 
\end{equation*}  
\begin{equation*}
+ H\, \Big( u\, \nabla^h_YZ + u\, H\, [vZ+wY] + [vw'+\alpha\alpha'\, h(Y,Z)]\, X \Big) . 
\end{equation*}  }\parn 
Therefore  
\begin{equation*}
\Big[\nabla^g_{(u\partial_t + X)} , \nabla^g_{(v\partial_t + Y)}\Big] (w\partial_t + Z) 
\end{equation*}
\parn\vbox{ 
\begin{equation*}
= \Big( [uv'-u'v] w' +(\alpha\alpha')'\, [u\, h(Y,Z) - v\, h(X,Z)] + \alpha\alpha'\,[X h(Y,Z) - Y h(X,Z)] \Big)\partial_t 
\end{equation*}
\begin{equation*}
+\,\Big(\alpha\alpha'\, [ h(X,\nabla^h_YZ) - h(Y,\nabla^h_XZ) ] + \alpha^2  H^2\, [v\,h(X,Z) - u\,h(Y,Z)] \Big) \partial_t 
\end{equation*}
\begin{equation*}
+ [\nabla^h_X , \nabla^h_Y]Z + H\, (v\,\nabla^h_XZ - u\,\nabla^h_YZ+w\,[X,Y]) + (Hw)'\, [u\,Y - v\, X] + H\,[uv'-u'v] \,Z 
\end{equation*} 
\begin{equation*}
+ H\, \Big( u\,\nabla^h_YZ - v\,\nabla^h_XZ + (Hw-w') (u\,Y-v\,X) + \alpha\alpha'\,[ h(Y,Z)\, X - h(X,Z)\, Y]\Big)  .  
\end{equation*}  }\parn 
And since 
\begin{equation*}
\nabla^g_{[(u\partial_t + X), (v\partial_t + Y)]} (w\partial_t + Z) = \nabla^g_{([uv'-u'v]\,\partial_t + [X,Y])} (w\partial_t + Z) 
\end{equation*}
\begin{equation*}
= \Big( [uv'-u'v] w' + \alpha\alpha'\, h([X,Y],Z) \Big)\partial_t + \nabla^h_{[X,Y]}Z + H\, ([uv'-u'v]\, Z + w\, [X,Y])\, ,  
\end{equation*}
we get 
\begin{equation*}
\Big(\Big[\nabla^g_{(u\partial_t + X)} , \nabla^g_{(v\partial_t + Y)}\Big] - \nabla^g_{[(u\partial_t + X), (v\partial_t + Y)]}\Big) (w\partial_t + Z) 
\end{equation*}
\vspace{-6mm} 
\parn\vbox{ 
\begin{equation*}
= \Big( (\alpha\alpha')'\, [u\, h(Y,Z) - v\, h(X,Z)] + \alpha\alpha'\,[X h(Y,Z) - Y h(X,Z)] \Big)\partial_t 
\end{equation*}
\begin{equation*}
+\, \Big( \alpha\alpha' [ h(X,\nabla^h_YZ) - h(Y,\nabla^h_XZ) -  h([X,Y],Z)] + \alpha^2  H^2\, [v\,h(X,Z) - u\,h(Y,Z)] \Big) \partial_t  
\end{equation*}
\begin{equation*}
+ \Big([\nabla^h_X , \nabla^h_Y]-\nabla^h_{[X,Y]}\Big)Z  + H\, (v\,\nabla^h_XZ - u\,\nabla^h_YZ) + (Hw)'\, [u\,Y - v\, X] \,  
\end{equation*} 
\begin{equation*}
 + H\, \nabla^h_{u\,Y-v\,X}Z + H\,(Hw-w') (u\,Y-v\,X) + \alpha^2H^2\,[ h(Y,Z)\, X - h(X,Z)\, Y] 
\end{equation*}  }\parn 
\begin{equation*}
=  \alpha \alpha'' h(uY- vX,Z) \partial_t + ([\nabla^h_X , \nabla^h_Y] - \nabla^h_{[X,Y]})Z + \frac{\alpha''}{\alpha} w (uY - v X) + \alpha'\!\,^2 [ h(Y,Z) X - h(X,Z) Y] .  
\end{equation*} 
By Formula (\ref{f.exprlocR}), this entails Formula (\ref{f.curvRwp}), which is equivalent to (\ref{f.curvRwploc}), by (\ref{f.exprlocR'}).\parsn 
Then, denoting by $(e_1,\ldots,e_d)$ an orthonormal basis of $(M,h)\,$: 
\parn\vbox{ 
\begin{equation*}
\langle {\rm Ricci}(v\partial_t + Y) \hbox{,}\, w\partial_t + Z\rangle_{\!\eta}\, = 
\Big\langle \RR\left(\partial_t\wedge (v\partial_t + Y)\right)\hbox{,}\, \partial_t\wedge (w\partial_t + Z)\Big\rangle_{\!\eta} \qquad 
\end{equation*} \vspace{-3mm} 
\begin{equation*}
\qquad -\, \alpha(t)\2 \sum_{j=1}^d \Big\langle \RR\left(e_j\wedge (v\partial_t + Y)\right)\hbox{,}\, e_j\wedge (w\partial_t + Z)\Big\rangle_{\!\eta}\, , 
\end{equation*} }\parn  
which yields (\ref{f.RicRwp})$\,$; which is in turn equivalent to (\ref{f.curvRicwploc}), by (\ref{f.exprlocRic}). $\;\diamond$ 

\bCor  \label{cor.wppf?} \   A Lorentz manifold $(\M,g)$ having the warped product form is of perfect fluid type (recall Section \ref{sec.perffluid}) \iff the Ricci operator of its Riemannian factor $M$ is conformal to the identical map$\,$: $\,{\rm Ric} = \Omega\times {\rm Id}\,$, \  for some $\,\Omega\in C^0(M)$. \parn 
If this holds, we must have$\,$: \  $U=\partial_t \,$,  $\,\tilde p + q = -d\,(\alpha''/\alpha)$,  and \  $q = \Omega\,\alpha\2 -(d-1) H' $.  
\eCor 
\ub{Proof} \quad  By (\ref{f.fpb'}) and (\ref{f.RicRwp}), this happens \iff for any $\,v\in C^0(I), Y\in \Gamma(TM)\,$: 
\begin{equation*} 
 \left\langle {\rm Ric}(Y),Y\right\rangle + [(d-1)\alpha'\,^2 + \alpha\alpha''] h(Y,Y) - d (\alpha''/\alpha) v^2 = q\, g(U,v\partial_t +Y)^2 + \tilde p \, [v^2-\alpha^2 h(Y,Y)] \, , 
\end{equation*} 
which (seen as a polynomial in $v$) forces $\,U=\partial_t \,$, and then splits into \  $\tilde p + q = -d\,(\alpha''/\alpha)$ and \  $\left\langle {\rm Ric}(Y),Y\right\rangle = -[(d-1)\alpha'\,^2 + \alpha\alpha'' + \tilde p\, \alpha^2]\,  h(Y,Y)$. This latter equation is equivalent to $\,  \tilde p = - [(d-1) H^2 + (\alpha'' /\alpha) + \Omega\,\alpha\2]\,$, and then, using $\,\tilde p + q = -d\,(\alpha''/\alpha)$, to $\,q = \Omega\,\alpha\2 -(d-1) H' $.  $\;\diamond$

\bCor  \label{cor.calcenwp} \  Consider a Lorentz manifold $(\M,g)$ having the warped product form. The energy (\ref{f.Ener}) at $\,\dot\xi\equiv (\dot t,\dot x)\in T^1_{\xi}\M \equiv  T^1_{(t,x)}(I\times M)$ equals$\,$: 
\begin{equation} \label{f.EnRwp} 
\EE(\xi,\dot\xi) = \langle {\rm Ric}(\dot x) , \dot x\rangle - (d-1)\, H'(t)\,({\dot t^2-1}) + \5\, d\, (d-1)\,  H^2(t) + {\frac{R^M}{2\, \alpha^2(t)}}\,\raise1.9pt\hbox{.} 
\end{equation}
Then the weak energy condition is equivalent to the following lower bounds for the Ricci operator and the scalar curvature $\,R^M$ of the Riemannian factor $(M,h)\,$: 
\begin{equation}  \label{f.EnRwpwec} 
\quad \inf_{w\in TM} \frac{ \langle {\rm Ric}(w) , w\rangle}{ h(w,w)} \ge  (d-1)\,\sup_{I}\,\{ \alpha\, \alpha'' - \alpha'\,^2\}  \; ;  \quad R^M\, \ge\, -\,d\, (d-1)\,  \inf_{I}\,\alpha'\,^2\, .
\end{equation}
And the scalar curvature $\,R\,$ of $(\M,g)$ equals$\,$: 
\begin{equation} \label{f.RicRwpf} 
R = - \alpha\2 R^{M} - d\Big[(d-1)\, |{\ts\frac{\alpha'}{\alpha}}|^2 + 2\, {\ts\frac{\alpha''}{\alpha}}\Big] ,  
\end{equation} 
so that its non-positivity is equivalent to the lower bound on the scalar curvature $\,R^M\,$: 
\begin{equation}  \label{f.SCRwpwec} 
R^M\, \ge -\, d\times  \inf_{I}\,\{ (d-1)\, \alpha'\,^2 + 2 \,\alpha\, \alpha''  \}\, .
\end{equation}
\eCor
\ub{Proof} \quad  From (\ref{f.RicRwp}), we compute the scalar curvature of $\M\,$: 
\begin{equation*} 
R = \langle {\rm Ricci}(\partial_t),\partial_t \rangle_{\!\eta} - \alpha\2 \sum_{j=1}^d 
\langle {\rm Ricci}(e_j), e_j\rangle_{\!\eta} 
= - \alpha\2 R^{M} - d\Big[(d-1)\, |{\ts\frac{\alpha'}{\alpha}}|^2 + 2\, {\ts\frac{\alpha''}{\alpha}}\Big] .  
\end{equation*} 
On the other hand, by (\ref{f.RicRwp}) we have at \   $\dot\xi\equiv (\dot t,\dot x)\in T^1_{\xi}\M \equiv  T^1_{(t,x)}(I\times M)\,$: 
\begin{equation*}
\langle {\rm Ricci}(\dot\xi), \dot\xi\rangle_{\!\eta} = \langle {\rm Ric}(\dot x), \dot x\rangle + [(d-1)\, |\alpha'|^2 + \alpha \alpha'']\, {\ts\frac{\dot t^2-1}{\alpha^2}} - d\, {\ts\frac{\alpha''}{\alpha}}\, \dot t^2 
\end{equation*}
\begin{equation*}
= \langle{\rm Ric}(\dot x),\dot x\rangle - (d-1)\, H'(t)\,({\dot t^2-1}) - d\,{\ts\frac{\alpha''}{\alpha}}(t)\, .
\end{equation*}
Thence Formula (\ref{f.EnRwp}). \   Then, the weak energy condition holds \iff for any $\,t\in I$ and $\,w\in TM\,$: 
\begin{equation*} 
\langle {\rm Ric}(w) , w\rangle \ge  (d-1)\, \alpha^2(t)H'(t)\, h(w,w) - \5\, d\, (d-1)\, H^2(t) - {\frac{R^M}{2\, \alpha^2(t)}}\,\raise1.9pt\hbox{.} 
\end{equation*}
By homogeneity with respect to $\,w\,$, this can be split into the following lower bound for the Ricci operator of the Riemannian factor $M\,$: 
\begin{equation*} 
\inf_{w\in TM}\,\frac{ \langle {\rm Ric}(w) , w\rangle}{ h(w,w)} \ge  (d-1)\,\sup_{I}\,\{ \alpha\, \alpha'' - \alpha'\,^2\}  \, ,  
\end{equation*}
together with the condition particularised to $\,w=0\,$, which yields the following lower bound for the scalar curvature of the Riemannian factor $M\,$: 
\begin{equation*} 
d (d-1)\, H^2(t) + \alpha(t)\2 R^M \, \ge 0\, , \; \hbox{ or }\quad  R^M\, \ge\, -\,d\, (d-1)\,  \inf_{I}\,\alpha'\,^2\, . \;\;\diamond 
\end{equation*}

\section{Example of Robertson-Walker (R-W) manifolds} \label{sec.exRW} \indf 
    These important manifolds are particular cases of warped product$\,$: they can be written $\,\M= I\times M\,$, where $\,I\,$ is an open interval of $\,\R_+\,$ and $\,M \in\{\S^3,\R^3,\H^3\}$, with 
spherical coordinates $\,\xi\equiv (t,r,\f, \psi)$ (which are global in the case of $\R^3,\H^3$, and are defined separately on two hemispheres in the case of $\S^3$), and are endowed with the pseudo-norm$\,$: 
\begin{equation}  \label{f.RWu} 
g(\dot\xi , \dot\xi) := \dot t^2 -\alpha(t)^2\bigg( \frac{\dot r^2 }{1-kr^2} + r^2 \dot \f^2 + r^2 \sin^2\!\f\,\, \dot \psi^2 \bigg) \, , 
\end{equation}  
where the constant scalar spatial curvature $\,k\,$ belongs to $\{-1, 0, 1\}$ (note that $\,r\in [0,1]$ for $\,k=1\,$ and $\,r\in \R_+$ for $\,k=0,-1$), and the expansion factor $\,\alpha\,$ is as in the previous section \ref{sec.CWP}.  \  Note that we have necessarily $\,\dot t\ge 1\,$ everywhere on $T^1\M\,$. \parn  
By (\ref{f.curvRwp}), we have the curvature operator given by$\,$: 
\begin{equation*}
\Big\langle \RR\left((u\partial_t + X)\wedge (v\partial_t + Y)\right) \hbox{,}\, (a\partial_t + A)\wedge (w\partial_t + Z)\Big\rangle_{\!\eta} \qquad 
\end{equation*} 
\begin{equation*} 
\qquad =\, \alpha\alpha''\, h(uY-vX,aZ-wA) - \alpha^2(\alpha'\,^2+k)\, [h(X,A)h(Y,Z) - h(X,Z)h(Y,A)]\,. 
\end{equation*}
By (\ref{f.RicRwp}), the Ricci tensor $(\!(\tilde R_{\ell m})\!)$ is diagonal, with diagonal entries$\,$: 
\begin{equation*}
\Big(\! -3\,\frac{\alpha''(t)}{\alpha(t)}\, \raise 1.2pt\hbox{,}\, \frac{A(t)}{1-kr^2} \, \raise 2pt\hbox{,}\,  A(t)\, r^2 , A(t)\, r^2 \sin^2\!\f\Big) , \quad \hbox{where }\; A(t) := \alpha(t)\, \alpha''(t) + 2\, \alpha'(t)^2 + 2k\, , 
\end{equation*}
and the scalar curvature is \quad $R= -6\,[\alpha(t)\, \alpha''(t) + \alpha'(t)^2 + k]\, \alpha(t)\2$. \parn
  The Einstein energy-momentum tensor \  $\tilde R_{\ell m} - \5\,{{R}}\, g_{\ell m} = \tilde T_{\ell m}$ \    is diagonal as well, with diagonal entries$\,$:  
\begin{equation*} 
\Big( 3\,\frac{\alpha'(t)^2+k}{\alpha(t)^2}\, \raise 1pt\hbox{,}\, \frac{- \tilde A(t)} {1-kr^2} \, \raise 2pt\hbox{,}\,  - \tilde A(t)\, r^2 ,\, - \tilde A(t)\, r^2 \sin^2\!\f\Big) ,  \; \hbox{ with }\; \tilde A(t) := 2 \alpha(t) \alpha''(t) + \alpha'(t)^2 +k \, . 
\end{equation*} 
Hence, we have 
\begin{equation*} 
\tilde T_{\ell m} - \alpha(t)\2 \tilde A(t)\,  g_{\ell m} = 2\, [k\, \alpha(t)\2 -H'(t)]\, 1_{\{ \ell = m = 0\}}\, .  
\end{equation*} 
Thus, in accordance with Corollary \ref{cor.wppf?}, we have here an example of perfect fluid$\,$: Equation (\ref{f.fp}) holds, \   with 
\begin{equation} \label{f.EinsteinRW}  
U_j \equiv \delta^0_{j} \; ,\quad -p(\xi) = k\,\alpha(t)\2 + 2H'(t) + 3H^2(t) \; , \quad q(\xi) = 2\, [k\, \alpha(t)\2 -H'(t)]\, , 
\end{equation} 
\begin{equation*} \label{f.EinsteinRW'}  
\tilde p(\xi) = {-2}\,[ 2 k\,\alpha(t)\2 + H'(t) + 3H^2(t)]/({d-1})\,  .  
\end{equation*} 
\if{ %
\begin{equation*} 
\nabla_kU_{\ell} = -U_j\Gamma^j_{k\ell} = - \Gamma^0_{k\ell}\, , \quad \dot\xi_s^k\, \nabla_kU_{\ell}(\xi_s)\, \dot\xi_s^\ell = H(t_s) (1-\dot t_s^2)\, . 
\end{equation*} }\fi %
Note that  
\begin{equation} \label{f.energRW}  
\AA_s = U_i(\xi_s)\dot\xi^i_s = \dot t_s \quad \hbox{ and } \quad \EE_s = 2\, [k\, \alpha(t_s)\2 -H'(t_s)]\,\dot t_s^2 - p(\xi_s)\, . 
\end{equation} 
By Corollary \ref{cor.calcenwp} (or by Remark \ref{rem.RTpf}$(ii)$ as well), the weak energy condition is equivalent to $\,$: \quad $\alpha'\,^2 +k\, \ge \,  (\alpha\, \alpha'')^+$. \parm

   We shall consider only eternal Robertson-Walker space-times, which have their future-directed half-geodesics complete. This amounts to \   $I=\R_+^*\,$, together with \parn
${\ds \int^\infty \frac{\alpha}{\sqrt{1+\alpha^2}} =  \infty }\,$. \ 
   In the case of the basic relativistic diffusion (solving Equation (\ref{f.reldiff}) in such Robertson-Walker model), we have in particular$\,$: 
\begin{equation} \label{f.brddt} 
d \dot t_s = \rr \sqrt{\dot t_s^2-1}\, dw_s + {\ts\frac{3 \rr^2}{2}}\, \dot t_s\, ds  - H(t_s) [\dot t_s^2-1]\, ds \, .  
\end{equation} 

\subsection{$\Xi$-relativistic diffusions in an Einstein-de Sitter-like manifold} \label{sec.RWR1} \indf 
   We consider henceforth  the particular case \  $\,I =\,]0,\infty[\,$, $\,k= 0\,$,   and \   $\alpha(t) = t^c\,$,  \parn with exponent $\,c>0\,$.  Note that such expansion functions $\,\alpha\,$ can be obtained by solving a proportionality relation between $\,p\,$ and $\,q$ (see [H-E] or [L-L]). \parn 
Thus \quad  $q = 2c\, t\2\,$, \   $p = (2-3c) c \, t\2$, \  $R= -6\,c\,(2c-1)\, t\2\,$, \  $\EE = c\,t\2\,(2\,\dot t^2 +3c-2)\,$. \parsn 
Note that the weak energy condition holds. The scalar curvature is non-positive \iff  $\,c\ge 1/2\,$, and the pressure $\,p\,$ is non-negative \iff $\,c\le 2/3\,$.  \par
  Note that the particular case $\,c=\frac{2}{3}\,$ corresponds to a vanishing pressure $\,p\,$, and is precisely known as that of \ub{Einstein-de Sitter universe} (see for example [H-E]). And the analysis of [L-L] shows up precisely both limiting cases $\,c=\frac{2}{3}\,$ and $\,c=\frac{1}{2}\,$. 

\subsubsection{Basic relativistic diffusion in an Einstein-de Sitter-like manifold} \label{sec.RWR1o} \indf 
   In order to compare with the other relativistic diffusions, we mention first for the basic relativistic diffusion (of Section \ref{sec.reldiff}), the stochastic differential equations satisfied by the main coordinates $\,\dot t_s\,$ and  $\,\dot r_s\,$, appearing in the 4-dimensional sub-diffusion $(t_s, \dot t_s, r_s, \dot r_s)$. By (\ref{f.brddt}), we have, for independent standard real Brownian motions $w,\tilde w\,$: 
\begin{equation} \label{f.brddtc} 
d \dot t_s\, = \,\rr \sqrt{\dot t_s^2-1}\, dw_s + {\ts\frac{3 \rr^2}{2}}\, \dot t_s\, ds  - \frac{c}{t_s}\, (\dot t_s^2-1)\, ds \, ;
\end{equation} 
\begin{equation} \label{f.brddr} 
d \dot r_s = \frac{\rr \,\dot t_s\,\dot r_s}{\sqrt{\dot t_s^2-1}}\, dw_s + \rr \sqrt{\frac{1}{t_s^{2c}} -\frac{\dot r_s^2}{\dot t_s^2-1}} \, d\tilde w_s + {\ts\frac{3 \rr^2}{2}}\, \dot r_s\, ds  + \left[\frac{\dot t_s^2-1}{t_s^{2c}} - \dot r_s^2\right]\! \frac{ds}{r_s}  - \frac{2c}{t_s}\,\dot t_s\,\dot r_s\,ds  \, . 
\end{equation} 
Almost surely (see [A]), $\,\lim_{s\to\infty}\limits \dot t_s = \infty\,$, and $\,x_s/r_s \sim \dot x_s/| \dot x_s|\,$ converges in $\S^2$. \par   
     
\subsubsection{$R$-diffusion in an Einstein-de Sitter-like  manifold} \label{sec.RWR10} \indf 
   With the above, Section \ref{sec.RLRD}  reads here, for the $\,R$-relativistic diffusion, when $\,c\ge 1/2\,$: 
\begin{equation} \label{f.RRDsdeRW} 
d\dot\xi_s = \rr\, dM_s + 9 c\,(2c-1) \rr^2\, t_s\2\, \dot\xi_s\, ds - \G_{ij}^\cdot(\xi_s)\,  \dot\xi_s^i\, \dot\xi_s^j\, ds\,  , 
\end{equation} 
with the quadratic covariation matrix of the martingale part $\,dM_s\,$ given by$\,$: 
\begin{equation*}
\rr\2\, [ d \dot\xi^k_s , d \dot\xi^\ell_s] = 6c\,(2c-1)\, [\dot\xi^k_s\, \dot\xi^\ell_s - g^{k \ell}(\xi_s)]\, t_s\2\, ds\, , \quad \hbox{for }\; 0\le k,\ell\le d\,.
\end{equation*} 
In particular, we have for independent standard real Brownian motions $w,\tilde w\,$: 
\begin{equation} \label{f.RRDdt} 
d\dot t_s\, =\, \frac{\rr}{t_s} \sqrt{6c\,(2c-1) (\dot t_s^2-1)}\, dw_s + \frac{9  \rr^2 c\,(2c-1)}{t_s^2}\, \dot t_s\, ds - \frac{c}{t_s} (\dot t_s^2-1)\, ds\,  ;  
\end{equation} 
\begin{equation} \label{f.RRDdr} 
d \dot r_s\, =\, \frac{\rr \sqrt{6c\,(2c-1)}}{t_s} \left[\frac{\dot t_s\,\dot r_s}{\sqrt{\dot t_s^2-1}}\, dw_s + \sqrt{\frac{1}{t_s^{2c}} -\frac{\dot r_s^2}{\dot t_s^2-1}} \, d\tilde w_s\right] 
\end{equation} 
\begin{equation*}
\qquad  + \, \frac{9  \rr^2 c\,(2c-1)}{t_s^2}\,  \dot r_s\, ds  + \left[\frac{\dot t_s^2-1}{t_s^{2c}} - \dot r_s^2\right]\! \frac{ds}{r_s}  - \frac{2c}{t_s}\,\dot t_s\,\dot r_s\,ds  \, . 
\end{equation*} 
As the scalar curvature $\,R_s= 6c\,(1-2c)/t_s^2\,$ vanishes asymptotically, we expect that almost surely the $R$-diffusion behaves eventually as a timelike geodesic, and in particular that \  $\lim_{s\to\infty}\limits \dot t_s = 1\,$.

\subsubsection{$\EE$-diffusion in an Einstein-de Sitter-like manifold} \label{sec.RWR11} \indf 
   Similarly, using (\ref{f.dxipEdiff'}), (\ref{f.EinsteinRW}), (\ref{f.energRW}), we have here \   $\EE\,\dot\xi - \tilde T\dot\xi = 2(0-H')(\dot t^2\dot\xi - \dot t\,U)\,$,  \   
so that Section \ref{sec.RicciRelD} reads here, for the $\,\EE$-diffusion$\,$: 
\begin{equation} \label{f.RRDsdeEW} 
d\dot\xi_s = \rr\, dM_s + {\ts \frac{3 \rr^2 c}{2}} \, t_s\2\,(2\,\dot t_s^2 +3c-2)\, \dot\xi_s\, ds +  2 \rr^2 c  \, t_s\2 (\dot t_s\,\dot\xi_s - U_s)\, \dot t_s\,ds - \G_{ij}^\cdot(\xi_s)\,  \dot\xi_s^i\, \dot\xi_s^j\, ds\,  , 
\end{equation} 
with the quadratic covariation matrix of the martingale part $\,dM_s\,$ given by$\,$: 
\begin{equation*}
\rr\2 \,[ d \dot\xi^k_s , d \dot\xi^\ell_s] = c\,[\dot\xi^k_s\, \dot\xi^\ell_s - g^{k \ell}(\xi_s)]\, (2\,\dot t_s^2 +3c-2)\, t_s\2\,ds\, , \quad \hbox{for }\; 0\le k,\ell\le d\,.
\end{equation*} 
In particular, we have for some standard real Brownian motion $\,w\,$: 
\begin{equation} \label{f.equtk=0}
d \dot t_s = \frac{\rr\sqrt{c}}{t_s} \sqrt{(2\, \dot t_s^2-2+3c)(\dot t_s^2-1)}\, dw_s + c \left[ 5 \rr^2 (\dot t_s^2-1+ {\ts\frac{9c}{10}})\, {\frac{\dot t_s}{t_s^2}} - \frac{\dot t_s^2-1}{{t_s}}\right]  ds \,\raise0.2pt\hbox{;}  
\end{equation}    
\begin{equation} \label{f.ERDdr}
d\dot r_s = \frac{\rr\sqrt{c}}{t_s} \sqrt{2\, \dot t_s^2-2+3c}\, \left[\frac{\dot t_s\,\dot r_s}{\sqrt{\dot t_s^2-1}}\, dw_s + \sqrt{\frac{1}{t_s^{2c}} -\frac{\dot r_s^2}{\dot t_s^2-1}} \, d\tilde w_s\right] 
\end{equation} 
\begin{equation*}
\qquad  + \, \rr^2 c\, (5\,\dot t_s^2 - 3 + {\ts\frac{9c}{2}})\,\frac{\dot r_s}{t_s^2}\, ds - \frac{2 c}{t_s}\,\dot t_s\, \dot r_s\, ds  + \left[\frac{\dot t_s^2-1}{t_s^{2c}} - \dot r_s^2\right] \frac{ds}{r_s}\, \raise1.9pt\hbox{.} 
\end{equation*} 

\brem \label{rem.compXidiff} \  Comparison of $\,\Xi$-diffusions in an E.-d.S.-like manifold \parn {\rm   
 Along the preceding sections \ref{sec.RWR1o},  \ref{sec.RWR10}, \ref{sec.RWR11},  we specified to an Einstein-de Sitter-like manifold the various $\,\Xi$-diffusions we considered successively in Sections \ref{sec.reldiff}, \ref{sec.RLRD}, \ref{sec.RicciRelD}. Restricting to the only equation relating to the hyperbolic angle $\,\AA_s=\dot t_s\,$, or in other words, to the simplest sub-diffusion $(t_s,\dot t_s)$, this yields Equations (\ref{f.brddtc}),  (\ref{f.RRDdt}),  (\ref{f.equtk=0}) respectively. We observe that even in this simple case, all these covariant relativistic diffusions differ notably, having pairwise distinct minimal sub-diffusions (with 3 non-proportional diffusion factors).
 }\erem 
\vspace{-3mm} %

\subsection{Asymptotic behavior of the $\,R$-diffusion in an \\ Einstein-de Sitter-like manifold} \label{sec.SasSt}   \indf 
   We present here the asymptotic study of the $\,R$-diffusion of an Einstein-de Sitter-like manifold (recall Sections \ref{sec.RWR1}, \ref{sec.RWR10}). We will focus our attention on the simplest sub-diffusion $(t_s,\dot t_s)$, and on the space component $\,x_s\in\R^3\,$. Recall from (\ref{f.energRW}) that $\,\dot t_s=\AA_s\,$ equals the hyperbolic angle, measuring the gap between the ambient fluid and the velocity of the diffusing particle. Recall also that, by the unit pseudo-norm relation, $\,\dot t_s\,$ controls the behavior of the whole velocity $\,\dot \xi_s\,$. We get as a consequence the asymptotic behavior of the energy $\,\EE_s\,$. \quad  
   As quoted in Section \ref{sec.RWR10}, we must have here $\,c\ge \5\,$\raise1.1pt\hbox{.} \parsn  
Note that for $\,c=\5\,$\raise0.7pt\hbox{,} the scalar curvature vanishes, and the $\,R$-diffusion reduces to the geodesic flow, whose equations are easily solved and whose time coordinate satisfies (for constants $\,a\,$ and $\,s_0$)$\,$: \parsn
\centerline{$s-s_0 = \sqrt{t_s\,(t_s+a^2)} - a^2 \log[\sqrt{t_s} + \sqrt{t_s+a^2}\,]\,$, \quad  whence \quad $t_s\sim s\,$.} \pars 
The proofs of this section (and of the following one) will use several times the elementary fact that \as\, a continuous local martingale cannot go to infinity. \pars 

   The following confirms a conjecture stated at the end of Section \ref{sec.RWR10}. 
\bpro \label{pro.dottto} \   The process $\,\dot t_s\,$ goes almost surely to 1, and $\,\EE_s\to 0\,$, as $\,s\to\infty\,$.  
\epro
\ub{Proof} \quad  By Equation (\ref{f.RRDdt}), 
\begin{equation*} 
\log\frac{\dot t_s}{\dot t_1} - 3\rr^2 c (2c-1)\! \int_1^s\! (2+\dot t_\tau\2)\frac{d\tau}{t_\tau^2} + c\! \int_1^s \! (1-\dot t_\tau\2)\, \frac{\dot t_\tau}{t_\tau}\, d\tau\, 
\end{equation*} 
is a continuous martingale with quadratic variation \  ${\ds 6\rr^2 c\, (2c-1)\! \int_1^s\! (1-\dot t_\tau\2) \frac{d\tau}{t_\tau^2}}\, $\raise1.8pt\hbox{.} 
Hence, since $\,\dot t_\tau\ge 1\,$ and therefore $\,t_\tau\ge \tau\,$, the non-negative process 
\begin{equation*} 
\log \dot t_s + c \int_1^s (1-\dot t_\tau\2)\, \frac{\dot t_\tau}{t_\tau}\, d\tau\, 
\end{equation*} 
converges \as\, as $\,s\to\infty\,$. 
This forces the almost sure convergence of the integral$\,$:    ${\ds \int_1^\infty (1-\dot t_\tau\2)\, \frac{\dot t_\tau}{t_\tau}\, d\tau < \infty\, }$, and of $\,\dot t_s\,$, towards some $\,\dot t_\infty\in\, [1,\infty[\,$. This implies in turn \  $ t_\tau = \O(\tau)$, hence \  ${\ds \int_1^\infty (1-\dot t_\tau\2)\, \frac{d\tau}{\tau}\, < \infty\, }$, whence finally $\,\dot t_\infty\,=1\,$.  $\;\diamond$      \pars 
 
     Consider now the functional \  $a := t^c\sqrt{\dot t^2-1}\,$, which is constant along any geodesic. 
\blem \label{lem.dotttoo0} \  For $\,c>\5$\raise0.7pt\hbox{,} the process \  $a_s := t_s^c\sqrt{\dot t_s^2-1}\,$ goes almost surely to infinity, and cannot vanish. Moreover, for any $\,\e>0\,$ we have almost surely$\,$: \  ${\ds \int_1^\infty  t_s^{2c-2}\,\frac{ds}{a_s^{2+\e}}\, < \infty}\,$. \parn
\elem 
\ub{Proof} \quad  We get from Equation (\ref{f.RRDdt})$\,$: 
\begin{equation*} 
da_s\, =\, \frac{\rr}{t_s} \sqrt{6c\,(2c-1)(a_s^2+ t_s^{2c})}\, dw_s + 3  \rr^2 c\,(2c-1)\,\frac{3\,a_s^2+ 2\,t_s^{2c}}{t_s^2\, a_s}\, ds \, ,  
\end{equation*} 
and then for any $\,\e\in\,]0,1]\,$ and for some continuous local martingale $\,M\,$: 
\begin{equation*}
0\le a_s^{-\e} = a_1^{-\e} - M_s - 3 \e \rr^2c\, (2c-1)\!\! \int_1^s \frac{[2-\e]\, a_\tau^2 + [1-\e]\, t_\tau^{2c}}{t_\tau^2\, a_\tau^{2+\e}}\, d\tau \raise0.1pt\hbox{.}   
\end{equation*}
The signs in this last formula, and the fact that \as\, a continuous local martingale cannot go to infinity, imply the convergence of the last integral and of the martingale term $\,M_s\,$,  entailing the almost sure existence of a finite limit $\,a_\infty\1\,$, hence of $\,a_\infty\in\,]0,\infty]$, and the almost sure convergence of the integral ${\ds \int_1^\infty \frac{d\tau}{a_\tau^{2+\e}\, t_\tau^{2-2c}}\,  < \infty}\,$. Now, by Proposition \ref{pro.dottto}, this implies ${\ds \int_1^\infty \frac{d\tau}{a_\tau^{2+\e}\, \tau} \le \int_1^\infty \frac{d\tau}{a_\tau^{2+\e}\, \tau^{2-2c}}\,  < \infty}\,$, hence $\,a_\infty =\infty\,$. \  Finally, the equation for $\,a_s^{-\e}\,$ forbids also the existence of a finite zero $\,s_0\,$ for $\,a_s\,$. Indeed, $\,s\nea s_0\,$ would force the martingale term of this equation to go to $-\infty$, which is impossible.  $\;\diamond$     
\parm

   The following reveals the  asymptotic behavior of the space component $(x_s)$ for $\,c>\5\,$\raise0.9pt\hbox{.} 
\bpro \label{pro.dotttOR} \   For $\,c>\5$\raise0.7pt\hbox{,} the space component converges \as\, (as $s\to\infty$)\,:   \parsn
\centerline{$\,x_s\, \to\, x_\infty\, \in \R^3$.} 
\epro 
\ub{Proof} \quad  $(i)$ \  Let us consider the non-negative process $\,u_s := t_s (\dot t_s^2-1)$, which is constant for $\,c=\5\,$\raise0.9pt\hbox{,} and cannot vanish for $\,c >\5\,$\raise0.9pt\hbox{,} by Lemma \ref{lem.dotttoo0}. By Equation (\ref{f.RRDdt}), 
\begin{equation*} 
u_s + (2c-1)\! \int_1^s \frac{u_\tau\,\dot t_\tau}{t_\tau}\, {d\tau}   
-\, 6\rr^2 c\, (2c-1)\! \int_1^s [4\, \dot t_\tau^2-1]\, \frac{d\tau}{t_\tau}\, 
\end{equation*} 
is a continuous local martingale. \  Then, for any  $\,\e >0\,$,  
\begin{equation*} 
\frac{u_s}{t_s^\e} - 6\rr^2 c\, [2c-1]\! \int_1^s\!  \frac{[4\dot t_\tau^2-1]\, d\tau}{t_\tau^{1+\e}} + [2c-1+\e]\!  \int_1^s\!\frac{u_\tau\, dt_\tau}{ t_\tau^{1+\e}}  
\end{equation*} 
is a continuous local martingale. By Proposition \ref{pro.dottto}, the central term converges \as. 
 This implies that  \  ${\ds  \int_1^\infty \frac{u_\tau\, dt_\tau}{ t_\tau^{1+\e}} < \infty}\;$ and that  $\,{\ds \frac{u_s}{t_s^\e}}\,$ converges, \as. \pars

   $(ii)$ \  By the unit pseudo-norm relation, we have \  $t_s^{2c}\, |\dot x_s|^2 = \dot t_s^2 - 1 = u_s/t_s\,$. Let us apply $(i)$ above with $\,\e:= c-\5\,$\raise0.7pt\hbox{,} to get\,: 
\begin{equation*}
\left[\int_1^\infty |\dot x_s|\, ds\right]^2 \le \int_1^\infty t_s^{\e-2c}\, ds \times \int_1^\infty t_s^{2c-\e}\,|\dot x_s|^2\, ds\, \le\, {\ts\frac{2}{2c-1}}  \int_1^\infty \frac{u_s}{t_s^{1+\e}}\, ds\, <\infty\, . 
\end{equation*}
This proves that \  ${\ds x_s = x_1 + \!\int_1^s\! \dot x_\tau\, d\tau \to  x_1 + \!\int_1^\infty\! \dot x_s\, ds \in \R^3}$, \as\, as $\,s\to\infty\,. \;\;\diamond$ \parm 

  In the case $\,c=\5\,$  of the $\,R$-diffusion being the geodesic flow, we have 
\begin{equation*}
r_s = \sqrt{b^2/a^2+ (a+o(1))\log s }\;\sim\, \sqrt{a \log s}\; \quad \hbox{ as }\; s\to\infty\, ,  
\end{equation*}
which shows that Proposition \ref{pro.dotttOR} does not hold for the limiting case $\,c=\5\,$\raise0.7pt\hbox{.} \pars
    To compare the $\,R$-diffusion with geodesics, note that (as is easily seen\,; see for example [A]) along any timelike geodesic, we have \  ${\ds x_s = x_1 + \frac{\dot x_1}{|\dot x_1|} \int_1^s \frac{a\, d\tau}{t_\tau^{2c}}\,}$ \  \big(and $\, \frac{\dot x_s}{|\dot x_s|}= \frac{\dot x_1}{|\dot x_1|}$\big), which converges precisely for $\,c>\5\,$\raise0.6pt\hbox{;}  \   and along any lightlike geodesic, we have \  ${\ds x_s = x_1 + \frac{\dot x_1}{|\dot x_1|} \int_{t_1}^{t_s} \frac{d\tau}{\tau^{c}}\,\sim\, V\times s^{\frac{1-c}{1+c}}\, }$  \  \big(and $\, \frac{\dot x_s}{|\dot x_s|}= \frac{\dot x_1}{|\dot x_1|}$\big), which converges only for $\,c>1\,$.  \pars
    On the other hand, for $c\le 1\,$, the behavior of the basic relativistic diffusion proves to satisfy (see [A])\,: \quad 
${\ds r_s\;\mathop{\sim}_{s\to\infty}\limits\,  \int_1^s \frac{a_\tau\, d\tau}{t_\tau^{2c}}\,\lra \infty\,}$    (exponentially fast, at least for  $\,c<1$). \parm 
   Hence, the $\,R$-diffusion behaves asymptotically more like a (timelike) geodesic than like the basic relativistic diffusion. However, owing to Lemma \ref{lem.dotttoo0}, the asymptotic behavior of the $\,R$-diffusion seems to be somehow intermediate between those of the geodesic flow and of the basic relativistic diffusion. \par  

\subsection{Asymptotic energy of the $\,\EE$-diffusion  in an E.-d.S. manifold} \label{sec.BehavtCED}   \indf 
    We consider here the case  of Section \ref{sec.RWR11}, dealing with the energy diffusion in an Einstein-de Sitter-like manifold, and more precisely, with its absolute-time minimal sub-diffusion $(t_s,\dot t_s)$ satisfying Equation (\ref{f.equtk=0}), and with the resulting random energy\,:  \pars
\centerline{$\EE_s = c\,t_s\2\,(2\,\dot t_s^2 +3c-2) = 2c\, (\dot t_s/t_s)^2 + \O(s\2)$.}\pars  
   Let us denote by  $\,\zeta\,$ the explosion time\,: \quad $\zeta := \sup\{ s>0\,|\, \dot t_s < \infty\}\in \,]0,\infty]\,$. \par 
\blem \label{pro.dottto1} \   We have almost surely$\,$:  \  either  \  $\lim_{s\to\zeta}\limits\, \dot t_s =1\;$  and $\;\zeta=\infty\,$, \  or \  $\lim_{s\to\zeta}\limits\, \dot t_s =\infty\,$. 
\elem
\ub{Proof} \quad  By Equation (\ref{f.equtk=0}), 
\begin{equation*}  
\frac{1}{\dot t_{s\wedge\zeta}} - c\!\int_{s_0}^{s\wedge\zeta}\! \big[1-{\ts\frac{1}{\dot t_\tau^2}}\big] \frac{d\tau}{{t_\tau}} + \rr^2 c \!\int_{s_0}^{s\wedge\zeta}\! \big[3+{\ts\frac{3c+2}{2\, \dot t_\tau^2}} + {\ts\frac{3c-2}{\dot t_\tau^4}}\big]  \frac{d t_\tau}{t_\tau^2} 
\end{equation*}  
is a continuous martingale with quadratic variation \  ${\ds 2\rr^2 c\! \int_{s_0}^{s\wedge\zeta}\! \big[1\!+\!{\ts\frac{3c-2}{2\, \dot t_\tau^2}}\big]\! \big[1\!-\!{\ts\frac{1}{\dot t_\tau^2}}\big] \frac{d\tau}{t_\tau^2}}\, $\raise1.8pt\hbox{.} 
Hence, since $\,\dot t_\tau\ge 1\,$ and therefore $\,t_\tau\ge \tau\,$, the process 
\begin{equation*} 
 \dot t\1_{s\wedge\zeta} - c \int_{s_0}^{s\wedge\zeta} (1-\dot t_\tau\2)\, \frac{d\tau}{t_\tau}\, 
\end{equation*} 
converges \as\, as $\,s\to\infty\,$. As $\,0\le \dot t\1_{s\wedge\zeta} \le 1\,$, this forces the almost sure convergence of the integral$\,$:    ${\ds \int_{s_0}^{\zeta} \big[1-{\ts\frac{1}{\dot t_\tau^2}}\big]\, \frac{d\tau}{t_\tau}\,  < \infty\, }$, and of $\,\dot t_{s\wedge\zeta}\,$, towards some $\,\dot t_{\zeta}\in\, [1,\infty]$.  Moreover, the convergence of the integral forces  either $\,\zeta<\infty\,$ and then $\,\dot t_{\zeta}=\infty\,$, or $\,\zeta=\infty\,$ and then 
$\,\dot t_{\zeta}\in\, \{1,\infty\}$. $\;\diamond$      \parm  
 
   The asymptotic behavior can, with positive probability, be partly opposite to that of the preceding $\,R$-diffusion$\,$: 
\bpro \label{pro.dottto2} \  From any starting point $(t_{s_0},\dot t_{s_0})$, there is a positive probability that both $\,\AA_s=\dot t_s\,$ and  the energy $\,\EE_s\,$  explode. This happens with arbitrary large probability, starting with $\,\dot t_{s_0}/t_{s_0}\,$ sufficiently large and $t_0$ bounded away from zero.  \parn
On the other hand, there is also a positive probability that the hyperbolic angle $\,\AA_s=\dot t_s\,$ does not explode and goes to 1, and then that the random energy $\,\EE_s\,$ goes to 0.  
This happens actually with arbitrary large probability, starting with sufficiently large $\,t_{s_0}/\dot t_{s_0}\,$. 
\epro
\ub{Proof} \quad  Let us set \  $\,\lambda_s := t_s/\dot t_s\, \ge 0\,$. From the above proof of Lemma \ref{pro.dottto1}, we get directly\,: 
\begin{equation} \label{dodo} 
\lambda_{s\wedge\zeta} - \lambda_{s_0} = M_{s\wedge\zeta}+ [1+c]\!\int_{s_0}^{s\wedge\zeta}\! \big[1-{\ts\frac{c}{[1+c]\dot t_\tau^2}}\big] {d\tau}  - \rr^2 c \!\int_{s_0}^{s\wedge\zeta}\! \big[3+{\ts\frac{3c-2}{2\, \dot t_\tau^2}} + {\ts\frac{3c-2}{\dot t_\tau^4}}\big]  \frac{d \tau}{\lambda_\tau} 
\end{equation}
\begin{equation} \label{dada}
\le \, M_{s\wedge\zeta} + [1+c] ({s\wedge\zeta}-s_0) - \rr^2 c \!\int_{s_0}^{s\wedge\zeta}\! \left[ \frac{3}{\lambda_\tau} - {\frac{1+2\,\dot t_\tau\2}{t_\tau\, \dot t_\tau}} \right]  {d \tau}\, , 
\end{equation}
$(M_s)$ denoting a martingale having quadratic variation \ 
${\ds 2\rr^2 c\! \int_{s_0}^{s}\! \big[1+{\ts\frac{3c-2}{2\, \dot t_\tau^2}}\big]\! \big[1-{\ts\frac{1}{\dot t_\tau^2}}\big] d\tau}$. 

   $(i)$ \  Let us first start the time sub-diffusion from $(t_{s_0},\dot t_{s_0})$ such that $\,t_{s_0}\ge {s_0}\ge 1\,$ and $\,\dot t_{s_0}\ge n\,m\,t_{s_0}\,$, with fixed $\, n\ge 2\,$ and $\, m\ge 2+ {\ts\frac{1+c}{3\rr^2 c}}\,$\raise1.5pt\hbox{,} \   and consider \parn    
$T:= \zeta\wedge\inf\{ s>s_0\,|\, m\, \lambda_s> 1\}$. \   Thus, we have on $[s_0,T]\,$: $\,\lambda_\tau\1\ge m\,$, and then \parn  
$\frac{3}{\lambda_\tau} - {\frac{1+2\,\dot t_\tau\2}{t_\tau\, \dot t_\tau}} \ge 3m-3 \ge 3+ {\ts\frac{1+c}{\rr^2 c}}\,$\raise1.5pt\hbox{.} \  Therefore, by  Inequality (\ref{dada}), we have \as \, for any $\,s\ge s_0\,$:  
\begin{equation*}  
0\le \, \lambda_{s\wedge T} \le  \lambda_{s_0} - 3 \rr^2 c\, ({s\wedge T}-s_0) + M_{s\wedge T} \, .   
\end{equation*}
Integrating this inequality and letting $\,s\nea \infty\,$ yields\,:  
\begin{equation*}
{\ts\frac{1}{m}}\,\P_{(t_{s_0},\,\dot t_{s_0})}[ T< \zeta ]\, \le \, \liminf_{s\to\infty}\, \E[\lambda_{s\wedge T} ] \,\le \lambda_{s_0} \le {\ts\frac{1}{nm}}\, \raise0.1pt\hbox{, \quad  whence }  \quad \P_{(t_{s_0},\,\dot t_{s_0})}[ T=\zeta ]\, \ge \, 1- {\ts\frac{1}{n}}\,\raise1.5pt\hbox{.}  
\end{equation*}
Moreover, \as\, on the event $\{T=\zeta\}$, the above inequality   
\begin{equation*}  
0\le \, \lambda_{s\wedge\zeta} \le  \lambda_{s_0} - 3 \rr^2 c\, ({s\wedge\zeta}-s_0) + M_{s\wedge\zeta} \,   
\end{equation*}
implies clearly (using that a continuous martingale \as\, cannot go to infinity) $\,\zeta < \infty\,$, and by the previous lemma that $\,\dot t_\zeta = \infty\,$. Then (\ref{dodo}) implies the convergence of $\,\lambda_{s\wedge\zeta}\,$ to some $\,\lambda_\zeta \in\R_+\,$.  \parn 
Furthermore, $\,\dot t_\zeta = \infty\,$ and $\,\lambda_\zeta >0\,$ for finite $\,\zeta\,$ would imply trivially $\,t_\zeta = \infty\,$, whence its logarithmic derivative should explode, which leads to a contradiction. \parn
This proves that we have \quad  ${\ds \P_{(t_{s_0},\,\dot t_{s_0})} \big[\zeta<\infty\,,\,\lambda_\zeta = 0\big]\, \ge \, 1-n\1 }$. \pars 
    Since (by the support theorem of Stroock and Varadhan, see for example Theorem 8.1 in [I-W]) from any starting point the sub-diffusion $(t_{s},\dot t_{s})$ hits with a positive probability some $(t_{s_0},\dot t_{s_0})$ as above, we find there is always a positive probability that 
$\,\dot t_s\,$ and $\,\EE_s\,$ explode (together). \pars 

   $(ii)$ \  Let us now start the time sub-diffusion from $(t_{s_1},\dot t_{s_1})$ such that  $\,n\, m'\,\dot t_{s_1} \le t_{s_1}\,$, with fixed $\, m'\ge 9\rr^2 c\,(1+c)\,$ and $\,n\ge 2+c\,$\hbox{,} \     and consider \     $T':= \zeta\wedge\inf\{ s>s_1\,|\, \lambda_s < m'\}$. \parn    
By Equation (\ref{f.equtk=0}) we have at once\,: \as, for any $\,s\ge s_1\,$, 
\begin{equation*} 
\lambda_{s\wedge T'}\1 = {\lambda_{s_1}\1} + 5\rr^2 c \!\int_{s_1}^{s\wedge T'}\! \big[1+{\ts\frac{9c-10}{10\, \dot t_\tau^2}}\big]  {\lambda_{\tau}^{-3}}\, d \tau
-[1+c] \!\int_{s_1}^{s\wedge T'}\! {\lambda_{\tau}\2}\, {d \tau} + c \!\int_{s_1}^{s\wedge T'}\! {t_{\tau}\2}\, {d \tau} + M'_{s\wedge T'}\, ,  
\end{equation*}    
$(M'_s)$ denoting a martingale having quadratic variation \ 
${\ds 2\rr^2 c\! \int_{s_1}^{s}\! \big[1+{\ts\frac{3c-2}{2\, \dot t_\tau^2}}\big]\! \big[1-{\ts\frac{1}{\dot t_\tau^2}}\big] \lambda_{\tau}^{-4}\, d\tau}$. \parn
Since on $[s_1,T']$ we have $\lambda_\tau\1\le 1/m'$, and then \  $5\rr^2 c \big[1+{\ts\frac{9c-10}{10\, \dot t_\tau^2}}\big] \lambda_\tau\1\le  \frac{5\rr^2 c\,(1+c)}{m'} \le 5/9\,$, \  we get\,:  
\begin{equation*} 
0 \le \, \lambda_{s\wedge T'}\1\, \le\, \lambda_{s_1}\1 - c\! \int_{s_1}^{s\wedge T'}\! {\lambda_{\tau}\2}\, {d \tau} + {\frac{c}{t_{s_1}}} + M'_{s\wedge T'}\, . 
\end{equation*}    
This entails \  ${\ds \int_{s_1}^{T'}\! {\lambda_{\tau}\2}\, {d \tau} <\infty\,}$ \  and  \  $\lambda_{s\wedge T'}\1 \to \lambda_{T'}\1 \in\, \R_+\,$ (which implies moreover \  $\,\lambda_{T'}\1 =0\,$ \as\, on $\{T'=\infty\}\,$), \  and 
\begin{equation*} 
{\ts\frac{1}{m'}}\, \P_{(t_{s_1},\dot t_{s_1})}[T'<\zeta] \le \E_{(t_{s_1},\dot t_{s_1})}\big[\lambda_{T'}\1\big] \le \lambda_{s_1}\1 + {\ts\frac{c}{t_{s_1}}} \le {\ts\frac{1+c}{n\,m'}}\, \raise1.7pt\hbox{.}  \end{equation*}    
Hence, we get \quad $ \P_{(t_{s_1},\dot t_{s_1})}[T'=\zeta]\, \ge 1-\frac{1+c}{n}\, >0\,$\raise0.4pt\hbox{.} \   
Furthermore, as in $(i)$ above, $\,\dot t_{T'} = \infty\,$ and $\,\lambda_{T'} >0\,$ for finite $\,T'\,$ is impossible, which excludes $\, T'=\zeta <\infty\,$. \   Therefore \parsn
\centerline{$ \P_{(t_{s_1},\dot t_{s_1})}[T'=\zeta=\infty]\, \ge 1-\frac{1+c}{n}\, >0\,$\raise0.4pt\hbox{.}}\parsn 
Then from the equation for $\,\lambda\,$, using that \   $\big[3+{\ts\frac{9c}{2}}\big] \rr^2 c\, \lambda_{\tau}\1\le \frac{(6+9c)\rr^2c}{2\,m'}< \5\, $ on $[s_1,T']$, we get \as\,:
\begin{equation*}  
\lambda_{s\wedge T'} - \lambda_{s_1} \ge ({s\wedge T'}-s_1) - \big[3+{\ts\frac{9c}{2}}\big] \rr^2 c \!\int_{s_1}^{s\wedge T'}\! \frac{d \tau}{\lambda_\tau}  + M_{s\wedge T'}\, 
\ge\,  \5\, (s\wedge T'-s_1) + M_{s\wedge T'}\, ,  
\end{equation*}  
which shows (since $[M_s,M_s] = \O(s)$) that \as\, $\{T'=\zeta =\infty\}\subset \{\lambda_s\to\infty\}$. \  On this same event, by Equation (\ref{f.equtk=0}) we have \as\, for any $\,s\ge s_1\,$:  
\begin{equation*} 
\dot t_s -\dot t_{s_1'} = \rr\! \int_{s_1'}^{s}\! \sqrt{2c \big[1+{\ts\frac{3c-2}{2\, \dot t_\tau^2}}\big]\! \big[1-{\ts\frac{1}{\dot t_\tau^2}}\big]}\, \frac{\dot t_\tau}{\lambda_{\tau}}\, d w_\tau
+ c \int_{s_1'}^{s} \left[ 5 \rr^2 \frac{1+ {\ts\frac{9c-10}{10\,\dot t_\tau^2}}}{\lambda_\tau^2} - \frac{1-\dot t_\tau\2}{{\lambda_\tau}}\right]  dt_\tau \,\raise0.2pt\hbox{,}  
\end{equation*}    
which shows that $\,\dot t_s\,$ cannot go to infinity, since this would forbid the last integral, and then the right hand side, to go to $ +\infty\,$. 
Hence, by Lemma \ref{pro.dottto1}, we obtain that 
\as\, $\{T'=\zeta =\infty\}\subset \{\dot t_s\to 1\}$. The proof is ended as in $(i)$ above, by applying the support theorem of Stroock and Varadhan, and by taking $\,n\,$ arbitrary large. $\;\diamond$

\section{Sectional relativistic diffusion}  \label{sec.SRD} \indf 
    We turn now our attention towards a different class of intrinsic relativistic generators on $\,G(\M)$,  whose expressions derive directly from the commutation relations of Section \ref{sec.FramebGM}, on canonical vector fields of $TG(\M)$. They all project on the unit tangent bundle $T^1\M$ onto a unique relativistic generator $\,\HH^1_{curv}\,$, whose expression involves the curvature tensor. Semi-ellipticity of $\,\HH^1_{curv}\,$ requires the assumption of non-negativity of timelike sectional curvatures. Note that in general $\,\HH^1_{curv}\,$ does not induce the geodesic flow in an empty space. 

\vspace{-3mm}  

\subsection{Intrinsic relativistic generators on $\,G(\M)$}  \label{sec.IntgenGM} \indf 
    We shall actually consider among these generators, those which are invariant under the action of $SO(d)$ on $G(\M)$.  \   To this aim, we introduce the following dual vertical vector fields, by lifting indexes$\,$: \  $V^{ij} := \eta^{im} \,\eta^{jn}\, V_{mn} \,$. \   Note that \  $\,V^j \equiv V^{0j} = - V_{0j} = - V_{j}\,$, and that \  $V^{ij} = V_{ij} \,$ for $\,1\le i,j\le d\,$. 
\    We consider again a positive parameter $\,\rr\,$.  
\bpro \label{pro.secreldiff} \   The following four $SO(d)$-invariant differential operators define the same operator $\,\HH^1_{curv}\,$ on $T^1\M\,$: 
\begin{equation*}  
H_0 - {\ts\frac{\rr^2}{2}} \sum_{j=1}^d \Big( [H_0,H_j] V^j +  V^j [H_0,H_j]\Big)\, ; \quad   H_0 + \rr^2\sum_{j=1}^d  [H_j,H_0] V^j \, ;\
\end{equation*}  
\begin{equation*}  
H_0 + \rr^2\sum_{j=1}^d R_{0}^{j}\, V_{j} - \rr^2 \sum_{1\le j,k\le d} \RR_{0}\!\,^{j0k}\, V_j V_{k} \; ;  \quad H_0 - {\ts\frac{\rr^2}{4}} \sum_{1\le i,j\le d} \Big( [H_i,H_j]\, V^{ij} + V^{ij}\, [H_i,H_j] \Big)\, ; 
\end{equation*}  

Note that $\,(\HH^1_{curv}- \LL^0)$ is self-adjoint with respect to the Liouville measure of $T^1\M$. 
\epro  
The proof will be broken in several lemmas. We begin with the following general and useful computation rules, derived from Sections \ref{sec.IsomL^2Rso1d} and \ref{sec.FramebGM}. 
\blem  \label{lem.calcViR} \quad   For $\,0\le j,k,\ell\le d\,$, \  we have$\,$: 
\begin{equation*}
V^i\RR_{ij}\!\,^{k\ell} = \delta_{0}^{k}\, R_j^{\ell} - \delta_{0}^{\ell}\, R_j^k +(1- d) \RR_{0j}\!\,^{k\ell} + \RR^{\ell}\!\,_{j0}\!\,^k - \RR^{k}\!\,_{j0}\!\,^\ell \, ; 
\end{equation*}  
\begin{equation*}
V^i\RR_{0i}\!\,^{k\ell} =  \delta_{0}^{\ell}\, R_0^k - \delta_{0}^{k}\, R_0^{\ell}\, ;  \qquad 
[\,[H_i,H_j], V^i]\,  = (d-1) [H_0,H_j]\, . 
\end{equation*}  
\elem 
\ub{Proof} \quad  We get the first formula by multiplying by $\,\eta^{ip}$ the formula of Lemma \ref{lem.calcViRo}, and particularising to $\,q=0\,$. As to the second one, by particularising the latter to $\,j=0\,$ and changing sign, we get :
\begin{equation*}
V^i\RR_{0i}\!\,^{k\ell} =  \delta_{0}^{\ell}\, R_0^k - \delta_{0}^{k}\, R_0^{\ell} + \RR^{k}\!\,_{00}\!\,^\ell - \RR^{\ell}\!\,_{00}\!\,^k . 
\end{equation*}  
Then, we note that \  $ \RR^{\ell}\!\,_{00}\!\,^k = \RR_{0}\!\,^{\ell k}\!\,_0 = \RR^{k}\!\,_{00}\!\,^{\ell}\!\,$.  \    Finally, the last formula derives from the second one and from the commutation relations (\ref{f.commrelVVH}) and (\ref{f.defRRTT}), as follows$\,$:
\begin{equation*}
[\,[H_i,H_j], V^i]\,  = \5\, [ \RR_{ij}\!\,^{k\ell}\, V_{k\ell} ,V^i] = \5\, \RR_{ij}\!\,^{k\ell}\, [V_{k\ell} , V^i] - \5\, (V^i\RR_{ij}\!\,^{k\ell})  V_{k\ell} 
\end{equation*}  
\begin{equation*}
=  \5\, \RR_{ij}\!\,^{k\ell} \Big(\delta_{k}^{i} V_{\ell} - \delta^{i}_{\ell} V_{k} + \eta^{ip}\, (\eta_{0\ell} V_{pk} - \eta_{0k} V_{p\ell})\Big)  
- \Big(\delta_{0}^{k}\, R_j^{\ell} + {\ts\frac{1-d}{2}}\, \RR_{0j}\!\,^{k\ell} + \RR^{\ell}\!_{j0}\!\,^k\Big) V_{k\ell}  
\end{equation*}  
\begin{equation*}
=  R_{j}^{\ell}\,V_{\ell} + \RR^p\!\,_{j}\!\,^{k}\!\,_{0}\, V_{pk} - R_{j}^{\ell}\,V_{\ell} + ({\ts\frac{d-1}{2}}) \RR_{0j}\!\,^{k\ell}\,V_{k\ell} - \RR^{\ell}\!_{j0}\!\,^k\, V_{k\ell}  
= ({\ts\frac{d-1}{2}}) \RR_{0j}\!\,^{k\ell}\, V_{k\ell} = (d-1) [H_0,H_j] . \;\diamond
\end{equation*}  
   
   We get then first the following.  
\blem \label{lem.F1} \   On  $\,C^2(T^1\M)$, we have  \   $ [\, [H_0,H_j], V^j] = 0\,$, \   and   
\begin{equation*}  
\HH^o_{curv} := - {\ts\frac{1}{2}}\sum_{j=1}^d \Big[ [H_0,H_j] V^j +  V^j [H_0,H_j]\Big] = \sum_{j=1}^d [H_j,H_0] V^j = \sum_{j=1}^d R_{0}^{j}\, V_{j} - \!\sum_{1\le j,k\le d}\!\! R_{0}\!\,^{j0k}\, V_j V_{k}\, .  
\end{equation*}  
\elem
\ub{Proof} \quad  Using the commutation relations (\ref{f.commrelVVH}) and (\ref{f.defRRTT}), that $\,V_{k\ell} = 0$ on $C^2(T^1\M)$ for $\,1\le k,\ell\le d\,$, and (\ref{f.RicciR}),  we have on one hand$\,$: 
\begin{equation*}
[H_0,H_j] V^j\, =\, \5\, \RR_{0j}\!\,^{k\ell}\, V_{k\ell}\, V^j\, = \,\5\, \RR_{0j}\!\,^{k\ell} ([V_{k\ell} , V^j] + V^j\,V_{k\ell})  
\end{equation*}  
\begin{equation*}  
= \5\, \RR_{0j}\!\,^{k\ell}\, \eta^{ij}\, (\eta_{ik} V_{\ell} - \eta_{i\ell} V_{k} + \eta_{0\ell} V_{ik} - \eta_{0k} V_{i\ell}) +  \RR_{0j}\!\,^{0\ell}\, V^j V_{\ell}\  
\end{equation*}  
\begin{equation*}  
=\, - \RR_{0j}\!\,^{kj}\, V_{k} + \RR_0\!\,^{ik}\!\,_0\, V_{ik} + \RR_{0j}\!\,^{0k}\, V^j V_{k}\, 
=\, \RR_{0j}\!\,^{0k}\, V^j V_{k} - R_{0}^{k}\, V_{k} \, . 
\end{equation*}  
On the other hand, using this first part of proof and Lemma \ref{lem.calcViR}, we get$\,$: 
\begin{equation*}
[ [H_0,H_j], V^j]\, =\, \5 [\RR_{0j}\!\,^{k\ell}\, V_{k\ell} , V^j]  \, 
=\, \5\, \RR_{0j}\!\,^{k\ell}\, [V_{k\ell} , V^j] - \5 (V^j\RR_{0j}\!\,^{k\ell})  V_{k\ell}\ 
\end{equation*}  
\begin{equation*}
=\, - R_{0}^{k}\, V_{k} - \5 (\delta_{0}^{\ell}\, R_0^k - \delta_{0}^{k}\, R_0^{\ell}) V_{k\ell}\, =\, - R_{0}^{k}\, V_{k} + R_{0}^{k}\, V_{k}\, = 0\, .  
\end{equation*}  
Using \  $\,[H_0,H_j] V^j +  V^j [H_0,H_j] = 2\,[H_0,H_j] V^j -  [\,[H_0,H_j] , V^j]\;$ ends the proof. $\;\diamond$
\parm   

   We get then the following. 
\blem \label{lem.F2} \  On  $\,C^2(T^1\M)$, we have  \   $\sum_{1\le i,j\le d}\limits [[H_i,H_j], V^{ij}] =0\,$,  \   and 
\begin{equation*}  
\sum_{1\le i,j\le d} \Big( [H_i,H_j]\, V^{ij} + V^{ij}\, [H_i,H_j]\Big) = \, -\, 4\, \HH^o_{curv}\; .   
\end{equation*} 
\elem
\ub{Proof} \quad As for the proof of Lemma \ref{lem.F1}, we use the commutation relations (\ref{f.commrelVVH}) and (\ref{f.defRRTT}), that $\,V_{k\ell} = 0$ on $C^2(T^1\M)\,$ for $\,1\le k,\ell\le d\,$, (\ref{f.RicciR}), Lemmas \ref{lem.calcViRo} and \ref{lem.calcViR}, and the symmetries of the Riemann tensor. \   We have thus on one hand and on $C^2(T^1\M)$$\,$:
\begin{equation*}  
[[H_p,H_j], V^{ij}] = \5 [\RR_{pj}\!\,^{k\ell} V_{k\ell} , V_{mn}]\, \eta^{im} \eta^{jn} = \5 \RR_{pj}\!\,^{k\ell} [V_{k\ell} , V_{mn}] \eta^{im} \eta^{jn} - \5 (V_{mn}\,\RR_{pj}\!\,^{k\ell}) \eta^{im} \eta^{jn} V_{k\ell}
\end{equation*} 
\parn \vbox{ 
\begin{equation*}  
=\, \5\, \RR_{pj}\!\,^{k\ell} \Big(\eta_{km}V_{\ell n} + \eta_{\ell n} V_{km} - \eta_{\ell m} V_{kn} - \eta_{kn} V_{\ell m}\Big) \eta^{im} \eta^{jn} \, -  \5 \,\times 
\end{equation*} 
\begin{equation*}
\Big[ \eta_{mp} \RR_{nj}\!^{k\ell} - \eta_{np} \RR_{mj}\!^{k\ell} + \eta_{mj} \RR_{np}\!^{k\ell} - \eta_{jn} \RR_{mp}\!^{k\ell} + \delta_{m}^{k} \RR_{pjn}\!^{\ell} - \delta_{n}^{k} \RR_{pjm}\!^{\ell} - \delta_{m}^{\ell} \RR_{pjn}\!^{k} + \delta_{n}^{\ell} \RR_{pjm}\!^{k} 
\Big]\! \eta^{im} \eta^{jn} V_{k\ell}
\end{equation*}  } 
\parn \vbox{ 
\begin{equation*}  
= \5\, \RR_{p}\!\,^{nk\ell} \Big(\delta_{k}^{i}\, V_{\ell n} + \eta_{\ell n}\, V_{k}\!\,^{i} - \delta^{i}_{\ell}\, V_{kn} - \eta_{kn} V_{\ell}\!\,^{i}\Big) 
\end{equation*} 
\begin{equation*}
- \, \5 \Big[-\RR^i\!\,_{p}\!\,^{k\ell} + \RR^i\!\,_{p}\!\,^{k\ell} - (d+1) \RR^i\!\,_{p}\!\,^{k\ell} + \eta^{ik} \RR_{pj}\!^{j\ell} - \eta^{jk} \RR_{pj}\!^{i\ell} - \eta^{i\ell} \RR_{pj}\!^{jk} + \eta^{j\ell} \RR_{pj}\!^{ik} \Big] \, V_{k\ell}
\end{equation*}  } \parn 
\begin{equation*}  
= \5 ( \RR_{p}\!^{ni\ell} V_{\ell n} + R_p^k\, V_{k}\!\,^{i} - \RR_{p}\!^{nki} V_{kn} + \RR_p^\ell\,  V_{\ell}\!\,^{i} ) + \5 [ (d+1) \RR^i\!\,_{p}\!\,^{k\ell} + \eta^{ik} R_{p}^{\ell} + \RR_{p}\!^{ki\ell} - \eta^{i\ell} R_{p}^{k} - \RR_{p}\!^{\ell ik}] V_{k\ell}
\end{equation*}  
\begin{equation*}  
=  \5 \Big[  - \RR_{p}\!\,^{ki\ell} + \eta^{i\ell} R_p^k - \RR_{p}\!\,^{\ell ki} - \eta^{ik} R_p^\ell + (d+1) \RR^i\!\,_{p}\!\,^{k\ell} + \eta^{ik} R_{p}^{\ell} + \RR_{p}\!^{ki\ell} - \eta^{i\ell}\, R_{p}^{k} - \RR_{p}\!^{\ell ik}\Big] V_{k\ell}
\end{equation*}  
\begin{equation*}  
=\,  ({\ts\frac{d+1}{2}})\, \RR^i\!\,_{p}\!\,^{k\ell} \, V_{k\ell}\, . 
\end{equation*}  
In particular, we get \quad 
$ [\,[H_i,H_j], V^{ij}] = ({\ts\frac{d+1}{2}})\, \RR^i\!\,_{i}\!\,^{k\ell} \, V_{k\ell} = 0\,$. \parsn 
And on the other hand, on $C^2(T^1\M)$ again$\,$:
\begin{equation*}  
[H_i,H_j] V^{ij}\, =\, \5\, \eta^{im} \eta^{jn}\, \RR_{ij}\!\,^{k\ell}\, V_{k\ell}\, V_{mn}\, =\, \eta^{i0} \eta^{jn}\, \RR_{ij}\!\,^{k\ell}\, V_{k\ell}\, V_{n}\, =\, \RR_{0}\!\,^{jk\ell}\, V_{k\ell}\, V_{j}  
\end{equation*} 
\begin{equation*}  
=\, \RR_{0}\!\,^{jk\ell} ([V_{k\ell} , V_{j}] +  V_{j}\, V_{k\ell})\,   
=\, \RR_{0}\!\,^{jk\ell} (\eta_{jk} V_{\ell} - \eta_{j\ell} V_{k} + \eta_{0\ell} V_{jk} - \eta_{0k} V_{j\ell})  + \RR_{0}\!\,^{jk\ell}\, V_{j}\, V_{k\ell}  
\end{equation*} 
\begin{equation*}  
= - 2\, R_0^k\, V_{k} + 2\, \RR_{0}\!\,^{jk}\!\,_0\, V_{jk} + 2\, \RR_{0}\!\,^{j0\ell}\, V_{j} V_{\ell}\,  
=\, 2\, (\RR_{0}\!\,^{j0k}\, V_{j} V_{k} - R_0^k\, V_{k})\, =\, - 2\, \HH^o_{curv} \, . \;\;\diamond 
\end{equation*} 

    The final assertion relating to the Liouville measure is proved as in Theorem \ref{the.gen}. 
\bpro \label{pro.F1}  \   In local coordinates, the second order operator $\,\HH^1_{curv}\,$ defined on $T^1\M$ by Proposition \ref{pro.secreldiff} reads\,: 
\begin{equation*}
\HH^1_{curv}\, =\, \dot\xi^j\frac{\p}{\p \xi^j} - \dot\xi^i\dot\xi^j\, \G_{ij}^k\, \frac{\p}{\p\dot\xi^k} +  {\ts\frac{\rr^2}{2}}\, \dot\xi^n \tilde R_{n}^{k}\, \frac{\p}{\p \dot\xi^k} - {\ts\frac{\rr^2}{2}}\,  \dot\xi^p\dot\xi^q\, \widetilde \RR_{p}\!\,^{k}\!\,_{q}\!\,^{\ell}\, \frac{\p^2 }{\p \dot\xi^k \p \dot\xi^\ell}
\end{equation*}
\begin{equation*}
=\, \dot\xi^j\frac{\p}{\p \xi^j} - \dot\xi^i\dot\xi^j\, \G_{ij}^k\, \frac{\p}{\p\dot\xi^k} +  {\ts\frac{\rr^2}{2}}\, \dot\xi^m\, \widetilde \RR_{mnpq}  \left( g^{nq} \, g^{pk}\, \frac{\p  }{\p \dot\xi^k} - \dot\xi^p\, g^{nk} g^{q\ell}\, \frac{\p^2 }{\p \dot\xi^k \p \dot\xi^\ell} \right)\raise1.8pt\hbox{.}
\end{equation*}
\epro 
\ub{Proof} \quad  By Section \ref{sec.ExprLocCoord}, we have on $C^2(T^1\M)$\,: 
\begin{equation*}
V_jV_k\, =\, e^n_je^\ell_k\,\frac{\p^2}{\p e^n_0\p e^\ell_0} + \d_{jk}\,e^n_0\,\frac{\p}{\p e^n_0} 
+ e^n_0e^\ell_0\,\frac{\p^2}{\p e^n_j\p e^\ell_k} + e^n_j\,\frac{\p}{\p e^n_k}\, =\, e^n_je^\ell_k\,\frac{\p^2}{\p e^n_0\p e^\ell_0} + \d_{jk}\,e^n_0\,\frac{\p}{\p e^n_0}\, \raise 2pt \hbox{,} 
\end{equation*}
whence by Lemma \ref{lem.F1}\,: 
\begin{equation*}
\HH^o_{curv}\, =\, - \RR_{0}\!\,^{j0k}\, e^n_je^\ell_k\,\frac{\p^2}{\p e^n_0\p e^\ell_0} - \RR_{0}\!\,^{j0k}\,\d_{jk}\,e^n_0\,\frac{\p }{\p e^n_0} + \sum_{j=1}^d\limits R_{0}^{j}\, e^n_j\,\frac{\p}{\p e^n_0} 
\end{equation*}
\begin{equation*}
\quad =\, - \RR_{0}\!\,^{j0k}\, e^n_je^q_k\,\frac{\p^2}{\p e^n_0\p e^q_0} + R_{0}^{j}\, e^n_j\,\frac{\p}{\p e^n_0} \quad \hbox{(including now $j=0$)}\, . 
\end{equation*}
On the other hand, by Formula (\ref{f.RtildeR}) we have$\,$: 
\begin{equation*}
\RR_{0}\!\,^{j0k} = \RR_{0abc} \,\eta^{aj} \eta^{b0} \eta^{ck} = e_0^m\, e_a^n\, \widetilde \RR_{mnp\ell}\, e^p_b\, e^\ell_c\,\eta^{aj}\, \eta^{b0}\, \eta^{ck} = e_0^m\, e^p_0\, e_a^r\, \widetilde \RR_{mrp\ell}\, e^\ell_c\,\eta^{aj} \, \eta^{ck} ,   
\end{equation*}
whence 
\begin{equation*}
\RR_{0}\!\,^{j0k}\, e^n_j\, e^q_k = e_0^m\, e^p_0\, e_a^r\, \widetilde \RR_{mrp\ell}\, e^\ell_c\,\eta^{aj} \, \eta^{ck}\, e^n_j\, e^q_k = e_0^m\, e^p_0\, \widetilde \RR_{mrp\ell}\, g^{rn} \, g^{q\ell} = e_0^m\, e^p_0\, \widetilde \RR_{m}\!\,^{n}\!\,_{p}\!\,^{q} .
\end{equation*}
And in a similar way, by (\ref{f.exprlocRic})$\,$:  
\begin{equation*}
R_{0}^{j}\, e^n_j = e_0^m\,e_i^q\,\tilde R_{mq}\,\eta^{ij}\, e^n_j = e_0^m\,\tilde R_{m}^{n} = e_0^m\, \widetilde \RR_{m\ell pq}\, g^{\ell q}\, g^{pn} .
\end{equation*}
This and (\ref{f.HV}),(\ref{f.V}) yield the wanted formula, whose coefficients depend only on $(\xi,\dot\xi)\in T^1\M$, as it must be by the $SO(d)$-invariance underlined in Proposition \ref{pro.secreldiff}. $\;\diamond $ \par

\subsection{Sign condition on timelike sectional curvatures}  \label{sec.PosCondCurv} \indf 
The generator $\,\HH^1_{curv}\,$ defined on $T^1\M$ by Proposition \ref{pro.secreldiff} is covariant with any Lorentz isometry of $(\M,g)$. \  Hence, it is a candidate to generate a covariant ``sectional'' relativistic diffusion on $T^1\M$, provided it be semi-elliptic.

   As a consequence of Section \ref{sec.IntgenGM}, the intrinsic sectional generator we are led to consider on $T^1\M\,$ is $\,\HH^1_{curv}\,$, restriction of  \  $H_0 + {\ts\frac{\rr^2}{2}} \big(R_{0}^{j}\, V_{j} - R_{0}\!\,^{j0k}\, V_j V_{k}\big)\,$. \parsn 

Now, a necessary and sufficient condition, in order that such an operator be the generator of a well-defined diffusion, is that it be subelliptic.   \parn  
   We are thus led to consider the following negativity condition on the curvature$\,$: 
\begin{equation}  \label{f.curvcond} 
\big\langle \RR(u\wedge v) , u\wedge v\big\rangle_\eta \, \le 0\, , \quad \hbox{ for any timelike $\,u\,$ and any spacelike $\,v\,$.} \; 
\end{equation} 
This condition is equivalent to the lower bound on sectional curvatures of timelike planes $\,\R u+\R v\,$: 
\begin{equation*}
\frac{\langle \RR(u\wedge v) , u\wedge v\rangle_\eta}{g(u\wedge v , u\wedge v)}\, \ge\, 0\, , 
\end{equation*}
since \   $\,g(u\wedge v , u\wedge v) := g(u,u)g(v,v)-g(u,v)^2 < 0$ \   for such planes. \parn
Note that Sectional Curvature has proved to be a natural tool in Lorentzian geometry, see for example [H], [H-R]. \par 
    We test this negativity condition on warped products, in Corollary \ref{cor.negcondwp} below.  \   When this negativity condition is fulfilled, we call the resulting covariant diffusion on $T^1\M$, which has generator $\,\HH^1_{curv}\,$ given by Propositions \ref{pro.secreldiff}, \ref{pro.F1}, the \ub{sectional relativistic diffusion}.  \par

\bcor  \label{cor.negcondwp} \  Consider a Lorentz manifold $(\M,g)$ having the warped product form. Then the sign condition (\ref{f.curvcond}) is equivalent to$\,$: \  ${\ds \alpha''\, \le 0}\,$ on $\,I$, \  together with the following lower bound on sectional curvatures of the Riemannian factor $(M,h)\,$: 
\begin{equation} \label{f.negcondwp} 
\inf_{X,Y \in TM} \, {\frac{\left\langle \KK\left(X\wedge Y\right) , X\wedge Y\right\rangle}{h(X,X)\, h(Y,Y)-h(X,Y)^2}}\,  \ge\, \sup_I\limits \{ \alpha\, \alpha'' - \alpha'\,^2\} .  
\end{equation} 
\ecor 
\ub{Proof} \quad   Let us denote by $\,SR(U\wedge V)$ the sectional curvature of the timelike plane 
associated with $\, U\wedge V$. Recall from Section \ref{sec.PosCondCurv} that the negativity condition (\ref{f.curvcond}) reads simply \   $SR(U\wedge V)\ge 0\,$, for any timelike $U$ and spacelike $V$. By choosing a pseudo-orthonormal basis of such given timelike plane, we can moreover restrict to $\,g(U,U) = 1 = - g(V,V)$ and $\,g(U,V)=0\,$. Setting $\,U= u \partial_t + X\,$ and  $\,V= v \partial_t + Y\,$, with $\,u,v\in C^0(I)$ and $\,X,Y\in TM$, we can thus suppose that$\,$: 
\begin{equation*}
u= \sqrt{\alpha^2\, h(X,X) + 1}\; ;\quad v= \sqrt{\alpha^2\, h(Y,Y) - 1}\; ;\quad uv= \alpha^2\, h(X,Y)\, ,
\end{equation*}
which implies \quad  $\alpha^2\, h(X\wedge Y,X\wedge Y) + h(Y,Y) - h(X,X) = \alpha\2$, \quad    and 
\begin{equation*}
h(uX-vY, uX-vY) = \alpha^2\, h(X\wedge Y,X\wedge Y) + \alpha\2\, . 
\end{equation*} 
Recall that $\,h(X\wedge Y,X\wedge Y) := h(X,X)h(Y,Y)-h(X,Y)^2$. \   Now, by (\ref{f.curvRwp}), this entails$\,$:   
\begin{equation*}
SR(U\wedge V) =  \alpha^2 \left\langle \KK\left(X\wedge Y\right) , X\wedge Y\right\rangle -  {\alpha \alpha''} [ \alpha^2\, h(X\wedge Y,X\wedge Y) + \alpha\2 ] + |\alpha \alpha'|^2\, h(X\wedge Y,X\wedge Y)  
\end{equation*} 
\begin{equation*}
= \alpha^2 \Big[ \left\langle \KK\left(X\wedge Y\right) , X\wedge Y\right\rangle - (\alpha\, \alpha'' - \alpha'\,^2) \, h(X\wedge Y,X\wedge Y)\Big] - \alpha''/\alpha \, . 
\end{equation*} \par 
   Reciprocally, given any $\,X,Y\in TM$ such that \parn   
   $\alpha^2\, h(X\wedge Y,X\wedge Y) + h(Y,Y) - h(X,X) = \alpha\2\,$  and \  $h(Y,Y) \ge \alpha\2\,$, \  setting \parn   
   $ u:= \sqrt{\alpha^2\, h(X,X) + 1}\;$ and  $\, v := \pm \sqrt{\alpha^2\, h(Y,Y) - 1}\,$, we get  \  $(uv)^2 = \alpha^4\, h(X,Y)^2\,$, \  and thence a basis $(U,V)$ as above. \parn 
Hence, the negativity condition (\ref{f.curvcond}) is equivalent to$\,$:  
\begin{equation*}
0 \le  \alpha^2 \Big( \left\langle \KK\left(X\wedge Y\right) , X\wedge Y\right\rangle - (\alpha\, \alpha'' - \alpha'\,^2) \, h(X\wedge Y,X\wedge Y) \Big) - \alpha''/\alpha\, , 
\end{equation*} 
for any $\,X,Y\in TM$ such that \  $\,\alpha^2\, h(X\wedge Y,X\wedge Y) + h(Y,Y) - h(X,X) = \alpha\2\,$  and \  $h(Y,Y) \ge \alpha\2\,$. \par 
   Distinguishing between collinear and non-collinear pairs $\,X,Y$, and denoting in the latter case by $\,SK(X\wedge Y)$ the sectional curvature of the plane associated with $\, X\wedge Y$, this condition splits into both$\,$: \quad $(\alpha''/\alpha) \le  0\,$, \quad together with$\,$: 
\begin{equation*}
\frac{(\alpha''/\alpha) }{\alpha^2\, h(X\wedge Y,X\wedge Y)}\; \le\,  SK(X\wedge Y) - (\alpha\, \alpha'' - \alpha'\,^2) \, , 
\end{equation*} 
for any non-collinear $\,X,Y\in TM$ such that \  $\,\alpha^2\, h(X\wedge Y,X\wedge Y) + h(Y,Y) - h(X,X) = \alpha\2\,$  and \  $h(Y,Y) \ge \alpha\2\,$. \   We shall have proved that this condition is equivalent to the wanted inequality (\ref{f.negcondwp}), if we show now that for any given $\,\alpha>0\,$, any given plane $P$ in $TM$ is generated by pairs $\,X,Y$ such that \  $\,\alpha^2\, h(X\wedge Y,X\wedge Y) + h(Y,Y) - h(X,X) = \alpha\2$, \  $h(Y,Y) \ge \alpha\2$, and such that $\,h(X\wedge Y,X\wedge Y)$ is arbitrary large. \parn
Now, starting from an arbitrary $\{X_0,Y_0\}$ generating $P$, \  take 
\  $Y := \frac{Y_0}{\alpha\sqrt{h(Y_0,Y_0)}}\,$ \  and \parn   
$X := q\Big[X_0 - \frac{h(X_0,Y_0)}{h(Y_0,Y_0)}\, Y_0\Big] $. Then \   
$\alpha^2\, h(X\wedge Y,X\wedge Y) + h(Y,Y) - h(X,X) = \alpha\2 = h(Y,Y)$, \  and    \   
$h(X\wedge Y,X\wedge Y) = q^2\, \frac{h(X_0\wedge Y_0,X_0\wedge Y_0)}{\alpha^2\, h(Y_0,Y_0)}\,$ is indeed arbitrary large, for arbitrary large $\,q\,$. $\;\diamond$ 
\pars

   In particular, in an Einstein-de Sitter-like manifold, the sign condition (\ref{f.curvcond}) holds  \iff  $\,\alpha''\le 0\,$, i.e. \iff $\,c\le 1\,$.
   
\vspace{-3mm}  

\subsection{Sectional diffusion in an Einstein-de Sitter-like manifold} \label{sec.RWR12} \indf 
   We must have here $\,c\le 1\,$. \   By (\ref{f.curvRwploc}), we have$\,$:   for $\,0\le k\le d\,$ and $\,1\le m,n,p,q\le d$,  
\begin{equation*} 
\widetilde \RR_{0nkq} =  \delta_{0k}\, c\,(c-1)\, t^{2c-2}\,h_{nq}\quad \hbox{ and } \quad 
\widetilde \RR_{mnpq} =  c^2\, t^{4c-2}\, [h_{mq}\, h_{np} - h_{mp}\, h_{nq}] - t^{2c}\, \tilde K_{mnpq}\, .  
\end{equation*}
Using Cartesian coordinates $(x^j)$ instead of spherical coordinates $(r,\f,\psi)$ for the Euclidian factor $\,M=\R^3$, we have merely$\,$: 
\begin{equation*} 
\widetilde \RR_{0nkq} = c\,(c-1)\, t^{2c-2}\, \delta_{0k}\, \delta_{nq}\quad \hbox{ and } \quad 
\widetilde \RR_{mnpq} =  c^2\, t^{4c-2}\, [\delta_{mq}\, \delta_{np} - \delta_{mp}\, \delta_{nq}] \, .  
\end{equation*}
Thence, \    for $\,0\le k\le d\,$ and $\,1\le m,n,p,q\le d\,$: 
\begin{equation*} 
\widetilde \RR_{0}\,\!^{n}\,\!_{k}\,\!^{q} = c(c-1)\, t^{-2c-2}\, \delta_{0k}\, \delta^{nq}\; \hbox{ and } \; \widetilde \RR_{m}\,\!^{n}\,\!_{p}\,\!^{q}  =  c(c-1)\, t^{2c-2}\,\delta_{mp}\,\delta^{tn}\, \delta^{tq} + c^2\, t^{-2}\, [\delta_{m}^{q}\, \delta_{p}^{n} - \delta_{mp}\, \delta^{nq}] .  
\end{equation*}
And \quad $ \tilde R_{ij} = \delta_{ij}\, \big(3c\,(1-c)\, t\2 \delta_{0j} + c\,(3c-1)\, t^{2c-2} [1-\delta_{0j}]  \big)$, \quad whence 
\begin{equation*} 
\tilde R_{n}^k = -\delta_{n}^{k}\,c\,t\2\, \big(3(c-1) \delta_{0n} + (3c-1) [1-\delta_{0n}]\big)\, . 
\end{equation*}
And by (\ref{f.LCcwploc}), the non-vanishing Christoffel coefficients are$\,$: 
\begin{equation*}  
\Gamma_{0j}^k = \Gamma_{j0}^k =  c\,t\1\, \delta_j^k\, , \quad \hbox{ and } \quad 
\Gamma_{ij}^0 = c\,t^{2c-1}\, \delta_{ij}\, , \quad \hbox{for }\; 1\le i,j,k\le d\,. 
\end{equation*}
Therefore, by Proposition \ref{pro.F1}, we have$\,$: 
\begin{equation*}
\HH_{curv}^1 = \dot\xi^j\frac{\p}{\p \xi^j} - \dot\xi^i\dot\xi^j\, \G_{ij}^k\, \frac{\p}{\p\dot\xi^k} +  {\ts\frac{\rr^2}{2}}\, \dot\xi^n \tilde R_{n}^{k}\, \frac{\p}{\p \dot\xi^k} - {\ts\frac{\rr^2}{2}}\,  \dot\xi^p\dot\xi^q\, \widetilde \RR_{p}\!\,^{k}\!\,_{q}\!\,^{\ell}\, \frac{\p^2}{\p \dot\xi^k \p \dot\xi^\ell}
\end{equation*}
\begin{equation*}
= \dot\xi^j\frac{\p}{\p \xi^j} - \frac{c}{t} (\dot t^2-1) \frac{\p}{\p\dot t} - \frac{2c}{t}\, \dot t\,\dot x^j \frac{\p}{\p\dot x^j} - \frac{3\rr^2c}{2\, t^2}(c-1)\, \dot t\, \frac{\p}{\p \dot t} - \frac{\rr^2c}{2\, t^2}(3c-1)\, \dot x^j \frac{\p}{\p \dot x^j}
\end{equation*}
\begin{equation*}
-\, {\frac{\rr^2c\,(c-1)}{2\,t^{2}}} (\dot t^2-1) \frac{\p^2}{\p\dot t^2} 
- {\frac{\rr^2 c\, (c-1)}{2\, t^{2c+2}}}\, \dot t^2\, \Delta_x 
+ {\frac{\rr^2c^2}{2\, t^2}}\,\left[ t^{-2c} (\dot t^2-1) \Delta_x - \dot x^i\dot x^j\, \frac{\p^2}{\p \dot x^i \p \dot x^j} \right] 
\end{equation*}
\begin{equation*}
=\,\dot\xi^j\frac{\p}{\p \xi^j} - \frac{c}{t} (\dot t^2-1) \frac{\p}{\p\dot t} - \frac{3\rr^2c}{2\, t^2}(c-1)\, \dot t\, \frac{\p}{\p \dot t} - \frac{2c}{t}\, \dot t\,\dot x^j \frac{\p}{\p\dot x^j} - \frac{\rr^2c\,(3c-1)}{2\, t^2} \, \dot x^j \frac{\p}{\p \dot x^j}
\end{equation*}
\begin{equation*}
+\,  {\frac{\rr^2c\,(1-c)}{2\,t^{2}}} (\dot t^2-1) \frac{\p^2}{\p\dot t^2} 
+ {\frac{\rr^2 c }{2\, t^{2c+2}}}\, (\dot t^2-c)\, \Delta_x 
- {\frac{\rr^2c^2}{2\, t^2}}\,\dot x^i\dot x^j\, \frac{\p^2}{\p \dot x^i \p \dot x^j} \, \raise2pt\hbox{.} 
\end{equation*}
We see that even in this simple case, the sectional and curvature diffusion differ tangibly, apart from the fact that the range of values of $c$ for which they are defined are different. 

\section{References} \label{sec.S} \pars  

\vbox{ \noindent 
{\bf [A]} \ Angst J.  \ \  {\it \'Etude de diffusions \`a valeurs dans des vari\'et\'es lorentziennes.}  
\par \hskip 19mm \quad Thesis of Strasbourg University, 2009. }
\par \medskip 

\vbox{ \noindent 
{\bf [A-C-T]} \  Arnaudon M., Coulibaly K.A., Thalmaier A. \  {\it Brownian motion with respect to\par \hskip 12mm a metric depending on time\,; definition, existence and applications to Ricci flow.}  \par \hskip 31mm	C. R. Acad. Sci. Paris  346, 773-778, 2008. } 
\parm 

\vbox{ \noindent 
{\bf [B]} \ Bailleul I.  \ \  {\it A stochastic approach to relativistic diffusions.}  
\par \hskip 21mm \quad Ann. I.H.P., vol. 46, n$^o\,$3, 760-795, 2010. }
\par \medskip 

\vbox{ \noindent 
{\bf [B-E]} \ Beem J.K.,  Ehrlich P.E. \ \  {\it Global Lorentzian Geometry.}
\par \hskip 56mm Marcel Dekker, New York and Basel, 1981. }
\par \medskip 

\vbox{ \noindent 
{\bf [De]} \ Debbasch F.  \ \  {\it A diffusion process in curved space-time.}  
\par \hskip 28mm \quad J. Math. Phys. 45, n$^o\,$7, 2744-2760, 2004. }
\par \medskip 

\vbox{ \noindent 
{\bf [Du]} \ Dudley R.M. \ \ {\it Lorentz-invariant Markov processes in relativistic phase space.}
\par \smallskip \hskip 33mm  Arkiv f\"or Matematik 6, n$^o\,$14, 241-268, 1965. }
\par \medskip 

\vbox{  \noindent 
{\bf [El]} \ Elworthy, D.  \  {\it Geometric aspects of diffusions on manifolds.} \par 
\par \hskip 16mm  {\'Ecole d'\'Et\'e de Probabilit\'es de Saint-Flour} (1985-87), Vol. XV-XVII, 277-425.  \par 
\par \hskip 16mm Lecture Notes in Math. Vol. 1362. Springer, Berlin, 1988.}
\parm 

\vbox{  \noindent 
{\bf [Em]} \ \'Emery, M.  \  {\it On two transfer principles in stochastic differential geometry.}
\par 
\par \hskip 16mm   {S\'eminaire de Pro\-babilit\'es} (1990), Vol. XXIV, 407-441. \  Springer, Berlin, 1990. }\parm 
\if{ %
\vbox{  \noindent 
{\bf [F1]} \ Franchi J.\quad {\it Asymptotic windings over the trefoil knot.}
\par \hskip 27mm   Revista Matem\'atica Iberoamericana, vol. 21, n$^o\, 3$, 729-770, 2005. }
\par \medskip  }\fi %

\vbox{  \noindent 
{\bf [F]} \ Franchi J.\quad {\it Relativistic Diffusion in G\"odel's Universe.}
\par \hskip 25mm   Comm. Math. Physics, vol. 290, n$^o\,$2, 523-555, 2009. }
\par \medskip  
 
\vbox{  \noindent 
{\bf [F-LJ]} \ Franchi J.,  Le Jan Y. \quad \hspace{-2mm}  {\it Relativistic Diffusions and Schwarzschild Geometry.} 
\par \hskip 53mm  Comm. Pure Appl. Math., vol. LX, n$^o\,$2, 187-251, 2007. }
\par \medskip  

\vbox{ \noindent 
{\bf [H-R]} \  Hall G.S.,  Rendall A.D. \quad  {\it Sectional Curvature in General Relativity.}
\par \hskip 57mm Gen. Rel. Grav., vol. 19, n$^o\,$8, 771-789, 1987. }
\par \medskip 

\vbox{ \noindent 
{\bf [H]} \  Harris S.G. \quad  {\it A Characterization of Robertson-Walker Spaces by Null Sectional\par  \hskip 29mm   Curvature.}
 \hskip 9mm Gen. Rel. Grav., vol. 17, n$^o\,$5, 493-498, 1985. }
\par \medskip 

\vbox{ \noindent 
{\bf [H-E]} \ Hawking S.W.,  Ellis G.F.R. \quad {\it The large-scale structure of space-time.}
\par \hskip 64mm Cambrige University Press, 1973. }
\par \medskip 

\vbox{ \noindent 
{\bf [Hs]} \ Hsu E.P.  \ \  {\it Stochastic analysis on manifolds. } \par \smallskip \hskip 25mm    
Graduate studies in Mathematics vol. 38, A.M.S., Providence 2002. }
\parm

\vbox{  \noindent 
{\bf [I-W]} \ Ikeda N.,  Watanabe S. \quad {\it Stochastic differential equations and diffusion
processes.}
\par \hskip 54mm North-Holland Kodansha, 1981. }
\par \medskip 

\vbox{ \noindent 
{\bf [K-N]} \ Kobayashi S.,  Nomizu K. \quad {\it Foundations of Differential Geometry.}
\par \hskip 9mm Interscience Publishers, John Wiley \& Sons, New York-London-Sydney, 1969. }
\par \medskip 

\vbox{  \noindent 
{\bf [L-L]} \ Landau L.,  Lifchitz E. \quad {\it Physique th\'eorique, tome II$\,$: Th\'eorie des champs.}
\par \hskip 52mm \'Editions MIR de Moscou, 1970. }
\par \medskip 

\vbox{  \noindent 
{\bf [M]} \ Malliavin, P. \   {\it G\'eom\'etrie diff\'erentielle stochastique.}  
\par \hskip 30mm    Les Presses de l'Universit\'e de Montr\'eal. \  Montr\'eal, 1978. }
\parm 

\centerline{--------------------------------------------------------------------------------------------}
\parn  
\ub{A.M.S. Classification} : \ Primary 58J65, secondary 53C50, 83C10, 60J65. \parsn 
\ub{Key Words} : \ Covariant relativistic diffusions, General relativity, Lorentz manifold, Riemann curvature tensor, Stochastic flow,  Robertson-Walker manifolds. 
\pars 
\centerline{--------------------------------------------------------------------------------------------}
\parmn 
\vbox{  \noindent 
$\bullet$ Jacques FRANCHI : \quad  jacques.franchi@math.unistra.fr \parn  
{\small \hskip0mm     IRMA, Universit\'e de Strasbourg et CNRS, 7 rue Ren\'e Descartes,  67084  Strasbourg cedex,  France\,; } \parn
b\'en\'eficiaire d'une aide de l'Agence Nationale de la Recherche, n$^o\,$ANR-09-BLAN-0364-01.  } 
\parm\parmn 
\vbox{  \noindent 
$\bullet$ Yves LE JAN : \quad yves.lejan@math.u-psud.fr \par 
\hskip7mm     Universit\'e Paris Sud 11, Math\'ematiques, B\^atiment 425, 91405 Orsay, France. }
\medskip 

\centerline{--------------------------------------------------------------------------------------------}

\end{document}